\documentclass[a4paper,12pt,reqno]{amsart}

\usepackage[mathscr]{euscript}
\usepackage{amssymb, amscd,enumerate,enumitem,color, mathtools, graphics, graphpap}
\usepackage[utf8x]{inputenc}
\usepackage[british]{babel}

\makeatletter
\def\subsection{\@startsection{subsection}{3}%
  \z@{.5\linespacing\@plus.7\linespacing}{.1\linespacing}%
  {\normalfont\bfseries}}
\makeatother

\font\smallsmc = cmcsc9
\font\smalltt = cmtt8
\font\smallit = cmti8

\numberwithin{equation}{section}

\theoremstyle{plain}
\newtheorem{theo}{Theorem}[section]
\newtheorem{lem}[theo]{Lemma}
\newtheorem{prop}[theo]{Proposition}

\newtheorem{cor}[theo]{Corollary}

\theoremstyle{definition}
\newtheorem{rem}[theo]{Remark}

\newtheorem{definition}[theo]{Definition}

\newenvironment{pf}{\noindent{\it Proof. }}{\hfill $\square$\par\medskip}
\newenvironment{pflemma}{\noindent{\it Proof of Lemma. }}{$\square$\par\medskip}

\theoremstyle{plain}

\theoremstyle{definition}

\newcommand{\beq}{\begin{equation}}
\newcommand{\eeq}{\end{equation}}
\renewcommand{\a}{\alpha}
\renewcommand{\b}{\beta}

\renewcommand{\d}{\delta}

\newcommand{\ve}{\varepsilon}
\newcommand{\f}{\varphi}
\newcommand{\g}{\gamma}

\renewcommand{\l}{\lambda}
\renewcommand{\o}{\omega}

\renewcommand{\r}{\rho}
\newcommand{\s}{\sigma}
\renewcommand{\t}{\tau}

\newcommand{\z}{\zeta}

\newcommand{\G}{\Gamma}

\renewcommand{\L}{\Lambda}
\renewcommand{\O}{\Omega}

\newcommand{\bC}{\mathbb{C}}
\newcommand{\bR}{\mathbb{R}}

\newcommand{\bH}{\mathbb{H}}

\newcommand{\bI}{\mathbb{I}}

\renewcommand{\gg}{\mathfrak{g}}
\newcommand{\gh}{\mathfrak{h}}
\newcommand{\gk}{\mathfrak{k}}

\newcommand{\gz}{\mathfrak{z}}

\newcommand{\gX}{\mathfrak{X}}

\newcommand{\gH}{\mathfrak{H}}

\newcommand{\su}{\mathfrak{su}}

\newcommand\GL{\mathrm{GL}}
\newcommand\SL{\mathrm{SL}}
\newcommand\SO{\mathrm{SO}}
\newcommand\SU{\mathrm{SU}}

\newcommand\Sp{\mathrm{Sp}}
\renewcommand\sp{\mathfrak{sp}}
\renewcommand\sl{\mathfrak{sl}}
\newcommand\ggl{\mathfrak{gl}}

\newcommand{\cA}{\mathscr{A}}
\newcommand{\cC}{\mathcal{C}}
\newcommand{\cD}{\mathscr{D}}
\newcommand{\cE}{\mathscr{E}}
\newcommand{\cF}{\mathscr{F}}
\newcommand{\cG}{\mathscr{G}}
\newcommand{\cH}{\mathscr{H}}

\newcommand{\cL}{\mathscr{L}}

\newcommand{\cO}{\mathscr{O}}
\newcommand{\cP}{\mathscr{P}}

\newcommand{\cS}{\mathscr{S}}
\newcommand{\cT}{\mathscr{T}}
\newcommand{\cU}{\mathscr{U}}
\newcommand{\cV}{\mathscr{V}}
\newcommand{\cW}{\mathscr{W}}

\newcommand{\Jst}{J_o{}} 

\newcommand{\p}{\partial}

\renewcommand{\square}{\kern1pt\vbox
{\hrule height 0.6pt\hbox{\vrule width 0.6pt\hskip 3pt
\vbox{\vskip 6pt}\hskip 3pt\vrule width 0.6pt}\hrule height0.6pt}\kern1pt}

\DeclareMathOperator\End{End\;}

\DeclareMathOperator\Ad{Ad}
\DeclareMathOperator\vol{vol}
\DeclareMathOperator\Id{Id}

\DeclareMathOperator{\Span}{span}
\DeclareMathOperator{\ad}{ad}
\DeclareMathOperator{\Diff}{Diff}

\renewcommand\={:=}

\newcommand{\wt}{\widetilde}
\newcommand{\wh}{\widehat}

\newcommand{\bt}{\begin{theo}\ \ }
\newcommand{\et}{\end{theo}}
\newcommand{\bp}{\begin{prop}\ \ }
\newcommand{\ep}{\end{prop}}
\newcommand{\bc}{\begin{cor}\ \ }
\newcommand{\ec}{\end{cor}}
\newcommand{\bl}{\begin{lem}\ \ }
\newcommand{\el}{\end{lem}}
\newcommand{\bd}{\begin{definition}}
\newcommand{\ed}{\end{definition}}
\newcommand{\n}{\nabla}

\newcommand{\be}{\begin{equation}}
\newcommand{\ee}{\end{equation}}
\newcommand\la[1]{\label{#1}}

\def\<#1,#2>{\langle\,#1,\,#2\,\rangle}
\newcommand{\arr}{\begin{array}{rlll}}
\newcommand{\ea}{\end{array}}
\newcommand{\bea}{\begin{eqnarray}}
\newcommand{\eea}{\end{eqnarray}}
\newcommand{\bean}{\begin{eqnarray*}}
\newcommand{\eean}{\end{eqnarray*}}

\newcommand{\HM}{{\cH(M)}}
\newcommand{\HCM}{{\cH^\bC(M)}}

\newcommand{\sH}{\mathbf H}
\newcommand{\sE}{\mathbf  E}

\def\sideremark#1{\ifvmode\leavevmode\fi\vadjust{
\vbox to0pt{\hbox to 0pt{\hskip\hsize\hskip1em
\vbox{\hsize3cm\tiny\raggedright\pretolerance10000
\noindent #1\hfill}\hss}\vbox to8pt{\vfil}\vss}}}
\begin{document}

\title[Instantons on hyperk\"ahler manifolds]
{Instantons 
on hyperk\"ahler manifolds}
\author[C. Devchand]{Chandrashekar Devchand}
\author[M. Pontecorvo]{Massimiliano Pontecorvo}
\author[A. Spiro]{Andrea  Spiro}

\begin{abstract} 
An instanton   $(E, D)$ on a (pseudo-)hyperk\"ahler manifold $M$  is a vector bundle $E$ associated to a  
principal $G$-bundle with a connection $D$ whose curvature is  pointwise invariant under the  quaternionic 
structures of  $T_x M, \ x\in M$, and thus satisfies the Yang-Mills equations. 
Revisiting   a construction of solutions, we prove a  local bijection between gauge equivalence classes of  instantons on $M$ 
 and equivalence classes of certain   holomorphic functions taking values in the Lie algebra of $G^\bC$ defined on an appropriate 
 $\SL_2(\bC)$-bundle over $M$.  Our reformulation affords a streamlined proof of  Uhlenbeck's   Compactness Theorem  for   instantons on (pseudo-)hyperk\"ahler manifolds.
\end{abstract}
\thanks{This research was partially supported by  
{\it Ministero dell'Istruzione, Universit\`a e Ri\-cer\-ca} in the framework of the project 
``Real and Complex Manifolds: Geometry, Topo\-logy and  Harmonic Analysis''  
and by GNSAGA of INdAM}

\maketitle
\section{Introduction}
\setcounter{equation}{0}
Beginning in the mid-1970's the self-duality equations for Yang-Mills fields successfully captured 
the imagination of theoretical physicists and mathematicians, epitomised by Donaldson's  flight into previously 
unforeseen realms of four-manifold differential topology (reviewed for instance in \cite{FU,Mo}). 
A Yang-Mills field is a pair  $(E, D)$ on a Riemannian manifold $(M, g)$, where  $E$  is a vector bundle 
associated to a  principle $G$-bundle with  a connection  $D$  whose curvature satisfies the Yang-Mills 
equation $D*F =0$. 
The (anti-)self-duality equations, requiring that the curvature $F$ of a connection $D$ over a Riemannian 
four-manifold $(M,g)$ takes values in the (anti-)self-dual eigenspace of the Hodge star-operator, 
implies the Yang-Mills equation  in virtue of the Bianchi identity $DF=0$.  
Connections satisfying the (anti-)self-dual Yang-Mills (SDYM) equations are called  (anti-) instantons. 
They are global minimisers of the Yang-Mills energy functional,
$S(A) = || F ||^2 = \int_M  F\wedge *F  \vol_g $.\par
 
The quest for explicit instanton solutions \cite{bpst} was initially physically motivated, for instance 
by the mystery of the phenomenon of quark confinement \cite{polyakov}, but the remarkable properties of
instantons soon attracted powerful mathematical treatment. First, Ward showed that solutions of the 
self-duality equations on $\bR^4$ are encoded in certain holomorphic  data on twistor space \cite{ward77},
effectively converting the problem to an algebro-geometric one. Then,  Atiyah, Hitchin and Singer \cite{ahs} obtained 
a correspondence between solutions of the SDYM equations on $S^4$ and certain real 
algebraic bundles on the complex projective 3-space $\bC P^3$. They thus established the relation between self-duality 
and holomorphic structures,  yielding in particular the dimension of the moduli space of solutions for any compact 
gauge group.  This led to a sequence of ans\"atze yielding SDYM solutions in terms of  arbitrary solutions 
of linear equations  \cite{aw, cfgy}.  Subsequently, powerful algebro-geometric results 
were used to obtain a complete construction of all SDYM fields on $S^4$ \cite{ADHM, DM, At}.   \par

These developments were followed by fundamental analytical results on variational methods for 
Yang-Mills theory. The moduli space of instantons is a subset of the quotient $\cA/\cG$ of the space of 
all connections $\cA$ with the group of all gauge transformations $\cG$. Locally representing the connections in   Coulomb gauges Uhlenbeck \cite{Uh1,Uh2} developed
analytical tools to study the singularities of the compact moduli space of instantons.   
Uhlenbeck's work, together with  the novel variational methods introduced by Taubes  to study gauge invariant theories, prepared the path for Donaldson's seminal work.  \par

These  analytical results depended crucially on the fact that  the Yang-Mills functional and therefore also the
Yang-Mills equations are conformally invariant in four dimensions. Further, the above-mentioned constructions of 
SDYM solutions crucially used the fact that $\bR^4$ conformally compactifies to  $S^4$. All this would seem to
impede any generalisation to Yang-Mills fields in higher dimensions. 
Indeed, it is known that a connection over the sphere  $S^d\,,\, d\ge 5$,  with sufficiently
small $L^2$-norm is necessarily flat \cite{BLS}.
However, the solvability of four dimensional 
SDYM  equations relies in particular on the fact that, being linear algebraic constraints on  
the curvature, they are first-order equations for the vector potential. This first-order property was partly
responsible for the good analytical properties of the SDYM equations. Indeed, an insistence upon this familiar sight of
partial-flatness conditions, requiring the vanishing of certain linear combinations of the curvature components, which 
automatically imply the second-order Yang-Mills equations,  yields the required instanton equations in dimensions 
greater than four. 
This idea was originally pursued in \cite{cdfn}, where it was shown that the required equations are restrictions of  
the curvature $F$ to an eigenspace of an endomorphism on the space of two-forms defined by an appropriate co-closed four-form $\O$,
\beq *(*\Omega \wedge F) = \l F \ ,\qquad  \O\in \L^4T^* M,\ \l\in \bR^*\ .
\la{gsd}\eeq 
The co-closedness of $\O$ suffices to show that a Yang-Mills  
curvature field satisfying \eqref{gsd} implies the second-order Yang-Mills equations.
For $d > 4$, the four-form $\O$ is pointwise  invariant only under
some proper subgroup  of  $ \SO_d(\bR)$. 
The existence of $\O$ corresponds to some special holonomy on the manifold \cite{DT, Do,UY, Ti, HL, ACD}.\par
 
The above-mentioned compactness results for the moduli space of Yang-Mills fields  were  already 
generalised to higher dimensions by Uhlenbeck and Nakajima \cite{Nak}.   For the  higher-dimensional 
generalisations of the self-duality equations the analytical programme in the spirit of Uhlenbeck and Taubes
was begun by Tian \cite{Ti}.  
The investigation of (local) solutions of higher dimensional equations of the form \eqref{gsd} began \cite{ward84,cgk,gios,ACD} 
with the case of instantons on spaces having hyperk\"ahler (hk) structure (i.e.\ with holonomy in $\Sp_n$), these being 
natural generalisations of $\bR^4 = \bH$. 
Some global results 
on general instantons on quaternionic K\"ahler (qk) manifolds (with  holonomy in  $ \Sp_n \cdot \Sp_1$) 
also exist; e.g. \cite{Ni,MS}. \par

The twistor formulation of SDYM, which led to the ADHM construction, has a (local) field theory variant, 
the {\it harmonic space formulation}. This was originally developed as a tool to study
supersymmetric harmonic maps \cite{gios, gios_book}. 
The harmonic space formulation enlarges the $\bC P^1$ fibre of the twistor bundle to $\SL_2(\bC)$, 
yielding a  total space more amenable to  standard field theoretical treatment. In the process 
(gauge equivalent classes of) local solutions of self-dual theories are parametrised by a prepotential, much as
the K\"ahler potential parametrises K\"ahler metrics. In a previous paper by two of us \cite{DS}, we have given 
a differential geometric description of the corresponding construction of (pseudo-)hyperk\"ahler metrics.  
In the current paper we investigate properties of  Yang-Mills instantons on (pseudo-)hyperk\"ahler manifolds 
using the harmonic space formulation, presenting a differential geometric formulation of the method based on 
the work of \cite{ACD}.  \par

The {\it harmonic space}  of an hk manifold $(M,g)$  is the trivial bundle  $\HM  = \SL_2(\bC) \times M \to M$,  
equipped with a  certain (non-product) complex structure. The space 
$\HM$ fibres naturally over the quotient   $Z(M) = \SL_2(\bC) /B\times M  \simeq \bC P^1 \times M$, 
where $B$ is  the Borel subgroup of upper triangular matrices. 
$Z(M)$ is the twistor bundle of  $(M, g)$ and has a  well defined complex structure, canonically determined by the hypercomplex structure of $M$. 
Now, the complex structure  of  the harmonic space $\HM$ is  the {\it unique} complex structure such that
 the  projection $p: \HM\to Z(M)$ is holomorphic. \par
 
A gauge field  $(E, D)$ on a complex hyperk\"ahler manifold $(M, g)$  is an instanton if the curvature of $D$  is  
pointwise invariant under the  quaternionic structure of  $T_pM, \ p\in M$.
Its pull-back field $(E', D')$ over $\HM$ admits an   {\it analytic gauge condition}, by which
we mean a special class of local trivialisations (= gauges) of $E'$. 
This class  has the crucial feature 
that its gauge  transformations are holomorphic, supplemented by some other conditions. 
In such a trivialisation, the potential   $A'$ of $D'$  is completely determined by just one of its 
components, which is moreover a holomorphic function on $\HM$.   
This component is called the {\it prepotential} of the gauge field $(E', D')$.
A freely-specifiable holomorphic prepotential, satisfying an appropriate first-order
linear equation on $\HM$, encodes all local properties of the corresponding instanton solution on $M$ and may be 
used to reconstruct the associated Yang-Mills field $(E, D)$. This construction, together with complete proofs of 
an essentially bijective correspondence between {\it normalised} prepotentials on $\HM$ and  moduli of  
locally defined instantons on  $M$ takes up the bulk of the content of this paper. 
\par

The existence of the analytic gauge condition and the resulting holomorphic prepotential
allows the use of  classical  results on holomorphic functions, such as Montel's Theorem or Hartogs' Removable Singularity Theorem,  to investigate the moduli spaces of instantons. 
Thus, in this formulation, holomorphy provides very useful tools. This is analogous to 
Uhlenbeck's Coulomb gauge condition, which allows the use of the machinery of elliptic equations.
As an example, we  establish some simple  estimates relating  $\cC^k$-norms   of prepotentials to those of  curvatures. These estimates,  combined  with the classical Montel Theorem  of   Complex Variable Theory,  
lead to a new direct proof of Uhlenbeck's  Strong Compactness Theorem for  instantons on hk manifolds   \cite{Uh1, Uh2, DK, Nak, Ti, We, Zh}.  \par
The paper is structured as follows. After the  preliminary section \S 2, we introduce  the notion of harmonic space and  discuss its relation with the twistor bundle of a (pseudo-)hyperk\"ahler manifold $M$ in \S 3. In \S 4, we discuss 
 the analytic gauge condition of the pull-back of the instanton over $\HM$ and the construction of the 
 instanton field over $M$ from the corresponding holomorphic prepotential on $\HM$. 
Our main new contributions appear in \S 4 and  \S 5, where we introduce a convenient normalisation for equivalent prepotentials,  prove a new existence result for an essentially unique instanton corresponding to a given prepotential,
obtain curvature estimates and present our new brief proof of the  Strong Compactness Theorem for  instantons on hk manifolds. 
\par

\noindent
{\bf Acknowledgements.}
One of us (CD) thanks Hermann Nicolai and the Albert Einstein Institute for providing an excellent research environment.
\par

\section{Preliminaries} 
\setcounter{equation}{0}

\subsection{Basics of  hyperk\"ahler manifolds} \label{adapted}
Given  a $4n$-dimensional real vector space $W$, we recall that a {\it hypercomplex structure\/}  on $W$ is a triple $(I_1, I_2, I_3)$ of endomorphisms  
satisfying the multiplicative relations of the imaginary quaternions,  
$I^2_\a = - \Id_W$ and $I_\a I_\b =  I_\g$
for all cyclic permutations  $(\a,\b,\g)$  of $(1,2,3)$.
Similarly, a  {\it hypercomplex structure} on a  $4n$-dimensional  real manifold $M$
is a triple $(J_1, J_2, J_3)$ of 
 integrable complex structures on $M$, with the property   that each triple $(I_\a \= J_\a|_x)$, $x \in M$, 
is a hypercomplex structure on  $T_xM$.\par
These two notions are   generalised as follows. Consider the $3$-dimensional subspace $Q_W$ of  $\End(W)$,
which is  the span  $Q_W = \Span_\bR(I_1, I_2, I_3)$ of the  endomorphisms of a hypercomplex structure  $(I_\a)$.
A {\it (pseudo-)quaternionic K\"ahler structure} on $W$
is an  inner product $\langle\cdot, \cdot \rangle$ on $W$  which is {\it hermitian with respect to $Q_W$},
that is, every element $ J \in Q_W$ is skew-symmetric with respect to $\langle\cdot, \cdot \rangle$.
For what concerns manifolds,  we have instead the following
\begin{definition}\label{defqk}
A $4n$-dimensional (pseudo-)Riemannian manifold $(M, g)$ 
of signature  $(4p, 4q)$, $p +q = n$,   is  a {\it  (pseudo-)quaternionic K\"ahler}  ({\it qk}) {\it manifold}
if it  admits a subbundle $Q \subset \End(TM)$ of quaternionic structures on the tangent spaces satisfying  the following two conditions: 
\begin{itemize} [itemsep=2pt, leftmargin=18pt]
\item[i)]   each inner product $g_x$ is Hermitian  for the quaternionic structure $Q_x$; 
\item[ii)]  the  parallel transport of the Levi-Civita connection $\n$ of $g$  preserves $Q$.
\end{itemize}
Further, $(M, g)$ is a {\it (pseudo-)hyperk\"ahler\/ (hk) manifold} if there exist three global $\nabla$-parallel sections 
$J_1, J_2, J_3$ of $Q$ (which are thus integrable complex structures) determining a hypercomplex structure on each tangent space of  $M$.
\end{definition}
It is well known that a qk manifold is Einstein. Moreover, it is hk if and only if its scalar curvature is  zero.
In this paper we  focus on hk manifolds, denoting them   by  a 5-tuple $(M, g, J_\a, \a = 1,2,3)$, where the $J_\a$  are  the
$\n$-parallel sections  generating the bundle $Q$.\par
 
 An  {\it adapted   frame} at the point $x\in M$ of an hk manifold $(M, g, J_\a)$ is a $g$-orthonormal frame 
 $u = (e_A): \bR^{4n}  \to T_xM$  mapping  the  standard triple of complex structures 
 $(\bf i, \bf j, \bf k)$  of $\bH^n \simeq \bR^{4n}$ into
 the  triple of complex structures $(J_\a|_x)_{\a = 1,2,3}$.  Note that   any  adapted frame  $u = (e_A)$ allows 
the identification of the  standard  $\Sp_1$- and $\Sp_{p,q}$-actions on  $\bH^n$ with  uniquely associated actions of $\Sp_1$- and $\Sp_{p,q}$-actions on $T_x M$.
The  family   $O_g(M, J_\a) $ of all adapted  frames is a principle bundle over $M$ with structure group  $\Sp_1{\cdot}\Sp_{p,q}$ and  is invariant under the parallel transport of the Levi-Civita connection.   
Hence, if  the bundle $O_g(M, J_\a)$ admits  a global section,  the above   
 actions of $\Sp_1{\cdot} \Sp_{p,q}$  on the  tangent spaces  $T_x M$  combine to  yield
 a {\it global} action of $\Sp_1{\cdot} \Sp_{p,q}$ onto $TM$.\par
\smallskip
We now fix   a few technical details, which we shall need. 
\par
As  the standard  representation of $\Sp_1{\cdot}\Sp_{p,q}$ on $\bH^n = \bC^{2n}$,  we choose the one  for   which 
the element  $ (\smallmatrix i & 0\\ 0 & -i  \endsmallmatrix )\in \sp_1$ acts on $(\bH^n)^\bC \simeq \bC^2 \otimes \bC^{2n}$ as the  left multiplication by  the complex structure $\mathbf i$. Hence,  for  any given  adapted frame $u \in O_g(M, J_\a)|_x$, the corresponding action of  $ (\smallmatrix i & 0\\ 0 & -i  \endsmallmatrix)$ 
 on $T_x M$  coincides with     the  action of the complex structure   $I = u_*(\mathbf i)$. 
This  action of  $\Sp_1 {\cdot} \Sp_{p,q}$ extends by $\bC$-linearity to  the
standard  action of  $\SL_2(\bC) {\cdot} \Sp_n(\bC)$ on
  $(\bH^n)^\bC \simeq \bC^2 \otimes \bC^{2n}$.  Consequently, the  $\bC$-linear extensions $u: (\bH^n)^\bC \to T^\bC_x M$ of  adapted frames give   
 $\SL_2(\bC) {\cdot} \Sp_n(\bC)$-equivariant isomorphisms $T^\bC_x M \simeq \sH_x \otimes \sE_x$  for 
 complex vector spaces     $\sH_x \simeq \bC^2$ and $\sE_x \simeq \bC^{2n}$. 
These isomorphisms are necessarily  related to each
 other by the  action  of  some element of  $ \Sp_1{\cdot}\Sp_{p,q} \subset \SL_2(\bC){\cdot} \Sp_n(\bC)$.\par
 Furthermore,  as  standard bases for $\bC^2$ and $\bC^{2n}$, we 
shall use 
$$\big(h^o_1 = (1,0), h^o_2 = (0,1)\big)\ ,\qquad \big(e^o_a = (0, \ldots, \underset{\text{$a$-th entry}}1, \ldots, 0)\big)\ .$$
For any ($\bC$-linearly extended) adapted frame $u:  (\bH^n)^\bC = \bC^2 \otimes \bC^{2n}  \to T^\bC_x M$, 
we then denote by $e_{ia}  \= u(h^o_i \otimes e^o_a)$, $i = 1, 2$, the complex vectors in $T^\bC_x M$ corresponding  to the
 basis 
$h^o_i \otimes e^o_a$ of $ \bC^2 \otimes \bC^{2n}$. 
If  $u: \bH^n \to T_x M$ is changed into   another adapted   frame 
 $$u' = u \circ U\qquad \text{with}\quad U = \left(\begin{array}{cc} u^1_+ & u^1_-\\ u^2_+ & u^2_-\end{array} \right)  \in \Sp_1\subset  \SL_2(\bC) = \Sp_1(\bC)\ ,$$
 the  corresponding complex basis $(e_{ja})$ of  $T^\bC_x M$ is changed  into 
 \beq \label{25} e_{+a} = u^j_+ e_{ja}\ ,\qquad e_{-a} = u^i_- e_{ia}\ .\eeq
Similarly,  if  $u' = u \circ U$ is further transformed  by   some  $A = (A_a^b) \in \Sp_{p,q} \subset \Sp_n(\bC)$, the  basis $(e_{\pm a})$ is then transformed into $e'_{\pm a} =A^b_a e_{\pm b}$. \par
 Note that, since the vectors  $h_1\otimes e_a^o$ and $h_2 \otimes e_a^o$ are $\pm i$-eigenvectors of the element $ (\smallmatrix i & 0\\ 0 & -i  \endsmallmatrix ) \simeq {\bf i} \in \Sp_1$, {\it the corresponding elements
   $e_{+a}$,  $e_{-a} \in T^\bC_x M$ are  vectors of type $(1,0)$ and $(0,1)$  with respect to the  complex structure $I = u(\bf i)$ of  $T_x M$}. \par
\subsection{Connections, gauges, potentials  and Yang-Mills fields}
In this section we briefly review certain basic  facts about  gauge theories 
which we shall need in what follows. 

Given a manifold $M$ and  a (principal or vector)  bundle $\pi: E \to M$, we denote by $\gX(M)$ the  space of all vector fields of $M$ and  by $\Gamma(E)$ the set of all sections $\s: M \to E$.  Given   $X \in \gX(M)$ and a smooth (real or complex) function $f$, we write  $X{\cdot} f$ to denote the directional derivative of $f$ along  $X$.\par
\smallskip
Let   $p: P \to M$  be a  principal $G$-bundle and  $p^E{:} E = P {\times_{G, \r}} V {\to} M$
 an associated vector bundle with fibre $V \simeq \bR^N$, determined by a faithful linear representation $\r: G \to \GL(V)$.
 We shall refer to an open subset $\cU \subset M$ on which there exists a choice of gauge (local trivialisation)   
 $\f: P|_{\cU} \to \cU \times G$ of the bundle $P$ over $\cU$ as the  {\it domain} of the   gauge $\f$ (a minor abuse of the language). 
Given  two gauges   $\f, \f'$ with  overlapping domains $\cU$, $\cU'$,     the transition function $\f' \circ \f^{-1}: (\cU \cap \cU') \times G \to 	(\cU\cap \cU') \times G$ is  a map of  the form 
$$\f' \circ \f^{-1}(x, h) = (x, g_x{\cdot} h)\qquad \text{for a smooth map }\quad g: \cU \to G\ .$$
The family of the  automorphisms of $G$ defined by    $h \mapsto g_x{\cdot} h$ is what is  usually called   {\it gauge transformation} between the gauges $\f$ and $\f'$. Finally, we recall  
that if a principal $G$-bundle   $p: P \to M$  admits a collection of gauges, whose domains
 form an open cover of $M$ and whose associated   gauge transformations 
 $h \mapsto g_x{\cdot} h$ are determined by maps $g: \cU \to G$ taking values in a fixed subgroup $G^o \subset G$, then such gauges determine  a $G^o$-subbundle $p^o: P^o \to M$ of $P$, called  {\it $G^o$-reduction}.  \par
\smallskip
A  connection $1$-form $\o$ on $P$ induces a unique covariant derivative  on the associated bundle $E$ . 
We recall that a {\it covariant derivative}  on  $E$ is an operator   
$D: \gX(M) \times \G(E) \to \G(E)$, which is   $\bR$-linear in both arguments,  $\cC^\infty(M; \bR)$-linear in 
$ \gX(M)$ and satisfies  the Leibniz rule
$D_X(f \s) = (X{\cdot} f) \s + f D_X \s$ for each  $f \in \cC^\infty(M; \bR)$.  

Now, consider  a gauge  $P|_{\cU}\overset \sim \to \cU \times G$ with domain $\cU$. 
In this gauge,   any vector  in 
$T P|_{\cU}$  is naturally identified with a sum $X_{(x,g)} +  B_{(x,g)} \in T_{(x, g)}( \cU \times G) \simeq T_x\cU + T_gG$ and the  $1$-form $\o$ on $P|_{\cU} \simeq \cU \times G$ can be pointwise expressed as a sum of  the form 
$$\o|_{(x, g)} = - A_{(x,g)} +  \o^G_g\ ,$$
where  $\o^G$ is the Maurer-Cartan form of   $G$ and, for each $g \in G$, the map   $A_{(\cdot, g)}: \cU \to T^* \cU \otimes \gg$
is a   $\gg$-valued $1$-form, which changes $G$-equivariantly with respect to $g$. 
The $1$-form  $A \= A_{(\cdot, e)}: \cU \to T^*M \otimes \gg$  is called   the
 {\it potential of $\o$ in the considered gauge}.  
 If   two   gauges   $P|_{\cU} \simeq \cU \times G$,   $P|_{\cU'} \simeq \cU' \times G$ have   overlapping domains $\cU, \cU'$, the corresponding     potentials $A$, $A'$  are related through   the gauge transformation $h \mapsto g_x{\cdot} h$  by means of 
 \beq \label{changepot} A' =  \Ad_{g^{-1}} A  + g^{-1}  {\cdot}  dg\ . \eeq
We now recall that a section $\s: \cU \to E|_{\cU}  \simeq \cU \times V $ of the associated bundle  has the form
$\s(x) = (x, s^i(x))$ for a smooth map $(s^i): \cU \to V = \bR^N$. Hence, 
for each  vector field  $X \in \gX( \cU)$, it is possible to  consider the   section  of $ E|_{\cU}  \simeq \cU \times V $
\beq \label{covder} D_X\s(x) =  \left(x, (X {\cdot} \s^i)|_x + A_j^i(X) \s^j|_x\right)\ , \eeq
 where $A^i_j \= \r_*{\circ}A$ with  $\r_*\colon \gg \to \ggl(V)$   
the Lie algebra representation determined by  the linear   representation $\r\colon G \to \GL(V)$ that gives the associated vector bundle $E = P \times_{G,\r} V$.  The  $\ggl(V)$-valued $1$-form $A^i_j$  is called the
 {\it  potential}  of $D$ induced by the connection $\o$. 
 The operator \eqref{covder} is gauge independent due to  \eqref{changepot}, thus   globally  well defined as a covariant derivative  on the sections of $E$.
 \par
\medskip
We recall that the  {\it curvature $2$-form of $\o$}  is the $\gg$-valued $2$-form  on $P$ defined by 
$ \O = d \o + \frac{1}{2}[\o, \o]$.  Given a gauge $\f: P|_{\cU} \to \cU \times \gg$ and the  associated  gauge $\wh \f$  for  $E$
$$ \wh \f: E|_{\cU} \to \cU \times V\ ,\qquad \wh \f\left([\f^{-1}(x, e), v)]\right) \= (x, v)\ ,$$  the  {\it curvature tensors} of  $\o$  and of  $D$   
 are  the ($\gg$- or $\ggl(V)$-valued)  $2$-forms $\cF^\f$ and $F_x$ on $\cU$,  defined   by 
\beq \label{2.4}
\cF^\f_x(v,w) \= 2 \O_{(x,e)}(v, w) \ ,\quad F^{\wh \f}_x(v, w) \=  2 \r_* \circ \cF^\f_x(v,w)\ , \quad\text{for }\  v, w \in T_x \cU\ .\eeq
Note that  $\cF^\f$ can be recovered from  the potential $A$  of $\o$  by the formula 
\beq \label{Fcurvature} \cF^\f(X, Y) = X {\cdot}(A(Y)) - Y {\cdot}(A(X)) + [A(X), A(Y)] - A([X, Y])\ .\eeq
and that,  if $h \mapsto g_x{\cdot} h$ is the  gauge transformation between  $\f$, $\f'$,  one has that 
\beq \label{transf1}  \cF^{\f'}_x = \Ad_{g_x} \cF^{\f}_x\ , \qquad F^{\wh \f'}_x = F^{\wh\f}_x\ . \eeq
This shows  that   {\it the curvature tensor  of} $\o$  {\it does depend  on the   gauge, while   the curvature  tensor of $D$  does not}.  Due to this, the curvatures  $F^{\wh\f}_x$  combine  and   determine a globally defined  $\ggl(V)$-valued  $2$-form $F$  on $M$. It can be also checked that it  
satisfies the identity 
\beq F(X, Y)s  = [D_X, D_Y]s - D_{[X,Y]} s\qquad \text{for sections}\ s \in \G(E)\ ,\eeq
a property often used as an alternative  definition of the curvature  of $D$.\par
\smallskip
All of the above notions and properties have analogues in  the 
case of holomorphic bundles, which we now   briefly recall. \par
If $(M, J)$ is a 
complex manifold,  $G$ is a complex Lie group and $V$ is a complex vector space, a principal $G$-bundle 
$p: P \to M$ (resp. a complex vector bundle $p: E \to M$ with fiber $V$) is called {\it holomorphic} if it is equipped with a complex structure $\wh J$, such that the right action of $G$ on $P$ (resp. the vector bundle structure on $E$) is $\wh J$-holomorphic  
and the projection $p$ is a $(\wh J, J)$-holomorphic mapping.  In this case, a trivialisation  is  called  {\it holomorphic} 
if it is a local holomorphic map from $(P, \wh J)$ (resp.  $(E, \wh J)$) to the cartesian product $M \times G$ (resp. $M \times V$), equipped with  the product complex structure.\par
A connection form $\o$ on a holomorphic $G$-bundle $(P, \wh J)$  is called {\it $\wh J$-invariant} if the corresponding horizontal spaces  $\cH_{(x,g)} = \ker \o_{(x,g)}$ are invariant under the complex structure $\wh J$. This is equivalent to say that,  in any holomorphic trivialisation $\f: P|_{\cU} \to \cU \times G$, the corresponding potential $A: \cU \to  T^* \cU \otimes  \gg$ takes actually values in in $ T^{10*} M \otimes \gg^{10} + T^{01*}  M \otimes \gg^{01}$. Here we denote by $ T^{10}_x M $ and $T^{01} M$ the holomorphic and anti-holomorphic distributions of $M$ and 
by $\gg^{10}$, $\gg^{01}$ the subalgebras of $\gg^\bC$ (both isomorphic to $\gg$), which are  generated by the vectors of type $(1,0)$ and $(0,1)$, respectively.  Being each $A_x$, $x \in \cU$,  real,  the projection of $A_x$  onto $T^{01*}_x M \otimes \gg^{01}$ is the complex conjugate of the component in $T^{10*}_x M \otimes \gg^{10}$. So,  the potential $A$ is uniquely determined by 
the associated  map   $A^{10}: \cU \to T^{10*} M \otimes \gg^{10}$, which we  call {\it $(1,0)$-potential}. \par
 We finally remark that the covariant derivative $D$ on an associated holomorphic vector bundle $p^E: E \to M$, determined by a $\wh J$-invariant  connection $\o$, is characterised by the property that it transforms sections of $E^\bC = E^{10} \oplus_M E^{01}$, with values in $ E^{10}$ or in $E^{01}$, into sections   which are still in  $ E^{10}$ or in $E^{01}$, respectively.  The covariant derivatives  of $(1,0)$-type (or $(0,1)$-type) sections and the change of 
  $(1,0)$-potentials  under holomorphic gauge transformations behave exactly as in  formulas  \eqref{covder}  and  \eqref{changepot}. \par  
\smallskip
We conclude this short section by recalling the definition of gauge fields. 
A {\it gauge field with structure group $G$} is a  pair $(E, D)$, formed by:
\begin{itemize}
\item[(1)] a vector bundle  $E$,  associated  with   a principal $G$-bundle $p: P \to M$,
\item[(2)] a covariant derivative $D$ on $E$, induced by a connection $\o$ on $P$. 
\end{itemize}
In this case we say that $(E, D)$ is a {\it   gauge field associated with  the pair $(P, \o)$}. If $G$ is complex,    $p:P \to M$ admits a reduction to a $G^o$- subbundle $P^o \subset P$  with    $G^o$ compact real form of $G$  and   $\o$ restricts  to a connection $\o^o$ on $P^o$, we  say that {\it $(P, \o)$ is reducible to $(P^o, \o^o)$} and that  {\it $(E, D)$ is  the {\rm complexification} of a gauge field with compact structure group $G^o$}.\par
\smallskip
 If $M$ is an oriented (pseudo-)Riemannian manifold, we may define the Hodge $\ast$-operator,  ${\ast}: \L^p T^* M \to \L^{n-p} T^* M$.
Then, the gauge field $(E, D)$ is called a {\it Yang-Mills field} if its curvature tensor $F$ satisfies the {\it Yang-Mills equation
$D {\ast} F = 0\,$}.\par
\subsection{$\cC^k$-norms and $L^p$-norms of curvatures}
\label{invariantnorms} 
 Let  $G$  be  a  reductive complex  Lie group and $G^o \subset G$  a compact real form of $G$.  Further let $\gg^o = \gg^s + \gz(\gg^o)$
  be the decomposition of $\gg^o = Lie(G^o)$  into its semisimple part  $\gg^s$  and center   $\gz(\gg^o)$,  and  denote by 
 $\langle \cdot, \cdot\rangle $  any    
  $\Ad_{G_o}$-invariant  Euclidean inner product  on  $\gg^o$ 
  which, on $\gg^s\times \gg^s$,   coincides with   minus   the Cartan-Killing form.  
This allows us to define the   Hermitian inner product  on
  $\gg = Lie(G) = \gg^o + i \gg^o$ 
\beq \label{innerproduct} h: \gg \times \gg \to \bC \ ,\quad h( X + i Y, X' + i Y' ) \= \langle X + i Y , X' - i Y' \rangle \eeq
and   the  associated norm   $\| X + i Y\|_h=\sqrt{\|X\|_{\langle \cdot, \cdot\rangle}^2 + \|Y\|_{\langle \cdot, \cdot\rangle}^2}$. 
If $(M, g)$ is    a Riemannian manifold, we can use the metric $g$  to extend  the Hermitian product   $h$  of $\gg$  to  a positive  inner product, also denoted by $h$, 
on the space of tensor fields in $\otimes^\ell T^* M\otimes \gg$  over $M$.  So,  for any compact subset $K \subset M$,  we  have 
the  usual    sup-norm  for  $\gg$-valued $\cC^k$ functions on $K$
$$\|V\|_{\cC^k(K, \gg)} \=  \sum_{j = 0}^k \sup_{x\in K}\|\n^{j} V|_x\|_{h}\ .$$
Similarly, for each   $p \in [1, +\infty)$, we may  use
the $L^p$-space $L^p(\cU, \gg)$,  the completion of the  space of  all $\cC^0$-maps  
 $V: \cU \to \gg$  with bounded values for the integral   
$ \int_\cU \|V\|^p_h\  \text{vol}_g$,  equipped with  the usual  $L^p$-norm  
$$\|\cdot \|_{L^p(\cU)} \=\left(\int_\cU \|\cdot\|_h^p\  \text{vol}_g\right)^{\frac{1}{p}}\ .$$
All such norms immediately generalise to  spaces of $\gg$-valued $r$-forms. \par
\smallskip
 Consider now a gauge field $(E, D)$  associated  with the pair $(P, \o)$ and a
trivialisation $\f: P|_{\cU} \to \cU \times G$ in which the  curvature tensors of 
$\o$ and $D$ are $\cF^\f$ and $F = F^{\wh \f}$, respectively. 
 Then, for any  compact subset $K \subset \cU$,  we  define  
 \beq \label{norms}\|F\|^{(\f)}_{\cC^k(K)} \= \|\cF^{\f}\|_{\cC^k(K, \gg)}\ ,\qquad  \|F\|^{(\f)}_{L^p(\cU)} \= \|\cF^{\f}\|_{L^p(\cU, \gg)}\ . \eeq
Now consider the cases when $(P, \o)$ is reducible to a pair $(P^o, \o^o)$ with structure group given by the compact $G^o$. Since the $h$-norms are $\Ad_{G_o}$-invariant, if we consider only the gauges which determine such a reduction,  then  the  norms \eqref{norms} do not depend on  $\f$  and  {\it  they coincide with  the usual $\cC^k$- and Sobolev norms  
of curvatures for  gauge fields with  compact structure groups}.   \par

\subsection{Lifting  gauge fields}
Let  $(E, D)$  be a gauge field on $M$  associated with a pair $(P, \o)$. 
If   $\pi\colon N {\to} M$  is  a principal $H$-bundle over  $M$, the {\it lift of $(P, \o)$}   is 
  the pair $(P',   \o')$ given by
 \begin{itemize}
 \item[a)] the lifted $G$-bundle  $p': P' \= \pi^* P \to N$,  i.e. the  submanifold  
$$\pi^* P\= \{(y, U)\in N \times P\ \text{such that}\  p(U) = \pi(y) \}\subset N \times P$$ 
equipped  with the natural projection 
$p': \pi^* P \to N$,   $p'(y, u) \= y$
\item[b)] the pull-back connection  $\o' \= \pi'{}^* \o$  on $P' = \pi^* P$ determined by  the projection  $\pi': P' \subset N \times P \to P$.
\end{itemize}
 The {\it lift of the gauge field $(E,D)$}  is 
the gauge field $(E' = \pi^* E, D' = \pi^* D)$ on $N$ given by 
\begin{itemize}
\item[i)] the lifted vector bundle $q':  E' = \pi^* E = P' \times_{G, \r} V \to N$ over $N$
\item[ii)] the covariant derivative $D' = \pi^* D$ on $E'$ induced  by $\o' $.
\end{itemize}
We now briefly discuss the problem of characterising the
gauge fields $(E', D')$ on principal $H$-bundles $\pi: N \to M$ which are  lifts of  gauge fields 
on  $M$. 
Given  $\pi: N \to M$  and $(E, D)$   as above,  for each  $X \in \gh = Lie(H)$,   the corresponding  $1$-parameter subgroup $\exp(\bR X) \subset H$ has  clearly a natural right action  on  $N$ and  determines a  natural $1$-parameter group of diffeomorphisms on the lifted bundle $P' \subset N \times P$ given by   
$$R': \bR \longrightarrow \Diff(P')\ ,\qquad R'_t(y, u) = (R_{\exp(t)}(y), u) \ .$$
 Note that each  map $R'_t$  is a    bundle automorphism that  preserves the connection $\o'$ and induces a  bundle automorphism of  the associated vector bundle $E' = \pi^* E$ that commutes  with  the covariant derivative $D' = \pi^* D$.
 Therefore, the  vector field $X'$  on $P'$, whose flow  is the  $1$-parameter group of automorphisms $R_t$,   is such that 
 \begin{itemize}
 \item[a)] $p'_*(X') = X$ and $\pi'_*(X') = 0$; here $X$ is identified with the   vertical vector field of the $H$-bundle $N$
 \item[b)] $\imath_{X'} \o' = 0$ 
and  $0 = \cL_{X'} \o' = \imath_{X'} d \o' + d(\imath_{X'} \o') = \imath_{X'} d \o'$; 
 \item[c)]  if    $A$ is the potential of $\o'$ in  a gauge  $\pi^* P|_{\cU} \simeq \cU \times G$,
 the  vector field $X'$ has  the form $X'_{(y, g)} = X_y + X^{\gg(y)}|_g$
  where $X^{\gg(y)}$ is the   left-invariant vector field   of $G$  defined by   
  $\o^G(X^{\gg(y)}) = A(X_y)$.
 \end{itemize}
Note that (b) implies also  that  
\beq  \label{2.8}  \imath_{X'}\O' = 0\qquad \text{and hence that }\qquad \imath_X F' =
 0\qquad  \text{for}\ X \in \gk\ .\eeq
 All this  has the following converse. \par
 \begin{prop} \label{prop28} Let  $(E', D')$ be  a gauge field associated with a pair $(P', \o')$ on   an $H$-bundle  $\pi: N \to M$ and, for any given gauge with domain $\cU \subset N$,  denote by  $\cF'$,  $F'$   the curvature tensors of $\o'$ and $D'$. \par
  If $H$ is simply connected, then 
   $(E', D')$ is the lift  of a gauge field on $M$ if and only if,  for  each $X \in \gh$,
    the associated infinitesimal transformation on $N$ is  such that
 a) $\imath_X \cF' = 0$ or, equivalently,  $\imath_X F' = 0$  in any gauge and  b) the unique $\o'$-horizontal vector field $X'$ on $P'$, which projects onto $X$,    is   complete.   \end{prop}
  \begin{pf} The necessity follows  from previous remarks. Assume now that for each $X \in \gh$, conditions (a) and (b)  hold. Since   $\imath_{X'} \o' = 0$, we have that 
  $\cL_{X'}  \o' = \imath_{X'} d \o' + d(\imath_{X'} \o') = \imath_{X'} \O' = 0$, 
  from which it follows  that the flow of $X'$  commutes with the right $G$-action of $P'$ and preserves $\o'$.  We also have that, for any  $X, Y \in \gh$, their $\o'$-horizontal lifts $X', Y'$ are such that 
  $\o'([X', Y']) =  \O'(X', Y') = 0$. This means that  $[X', Y']$ is the $\o'$-horizontal lift of  the Lie bracket $[X, Y]$.
    We therefore conclude that the collection of vector fields 
    $$\gh' \= \{\ X'\ \text{is $\o'$-horizontal  lift of some}\ \  X \in \gh\ \}$$
      is a finite dimensional Lie algebra of complete vector fields on $P'$.  By a classical theorem of Palais \cite{Pa},  this implies the existence of 
     a  right $H$-action  on $P'$   whose  infinitesimal transformations  are precisely the  $\o'$-horizontal vector fields $X' \in \gh'$.  All orbits of this action
      are regular and simply connected, because all of them  are  coverings of the simply connected $H$-orbits on $N$. Moreover,  each transformation of this action is   a bundle automorphism, which  preserves the connection $1$-form   $\o'$ and commutes with the $G$-action.    Hence,  the orbit space $P = P'/H$ is  a $G$-bundle  over $M \= P'/H' \times G = P/G$
 and is  equipped with the $\gg$-valued $1$-form $\o$ defined  by 
  $$\o_{[u]}(v) \= \o'_u(v')\  \text{for some}\ v' \in T_u P'\ \text{that projects onto}\  v\in T_{[u]} P'/H'\ .$$
  One can directly check that    $\o$ is a connection and that $\o' $ is the pull-back of $\o$ on $P'$. The associated  bundle $E$ of $P\to M = P/G$, equipped with the covariant derivation $D$ determined by  $\o$, is the desired gauge field,  of which $(E', D')$ is a  lift.
  \end{pf}
\section{Harmonic spaces of   (pseudo-)hyperk\"ahler manifolds} 
\setcounter{equation}{0}
In  the rest of this paper $(M, g, J_\a)$ denotes  an hk manifold and    $(E, D)$  a complex gauge field  on $(M, g, J_\a)$  associated with a pair $(P, \o)$  with complex structure group $G$.  We shall mostly assume that  $G$ is the complexification of a compact real form $G^o$ and that $(P, \o)$ is the complexification of a pair $(P^o, \o^o)$, with structure group $G^o$. We shall simply say  that $(E, D)$ is the ``complexification'' of a gauge field having compact structure   group $G^o$.  \par
\subsection{The twistor bundle of a (pseudo-)hyperk\"ahler manifold} \label{section321}
Let $(M, g, J_\a)$ be a $4n$-dimensional  hk manifold. It is well known that  for  a point  $z = (a, b, c) \in S^2  (\simeq \bC P^1)$, the tensor field $ I^{(z)} \= a J_1 + b J_2 + c J_3$ 
is an integrable complex structure on $M$ and that the {\it twistor bundle}  $Z(M)$ of $M$ is simply  the trivial $\bC P^1$-bundle  
$
Z(M) \=  M \times  \{I^{(z)}, z\in S^2\}  \to  M$.  Since each  $z \in S^2$ corresponds 
to a distinct integrable complex structure,  the twistor bundle $Z(M)$ is foliated by the complex submanifolds
$M \times\{z\} \simeq M$, each equipped with   the corresponding  complex structure  $I^{(z)}$. Such complex structures combine with the classical complex structure of $\bC P^1$  and determine  a natural almost complex structure on $Z(M)$, which we denote by $\wh I$. It was proved to be integrable in \cite{Sa}. \par
\smallskip
We remark  that the  complex structures on $T_x M,\ x \in M $, 
$$I^{(z)}|_x  \= a J_1|_x + b J_2|_x + c J_3|_x\ ,\qquad  z = (a, b, c) \in S^2 \simeq \bC P^1\ ,$$
coincide with the complex structures of the form $I = u_*(\mathbf i) \in \Span(I_\a \= J_\a|_x)$,  given  by  adapted frames $u$ of $T_x M$ as mentioned in  \S \ref{adapted}.  It follows that  each of the  $2n$-tuples of  complex vectors $e_{+a}$ (resp. $e_{-a}$), which are part of  the complex bases   \eqref{25},  is  actually  a frame  of holomorphic (resp. anti-holomorphic)  vectors 
for a complex leaf $M \times \{I^{(z)}\} \subset Z(M)$. 
  \par
\subsection{The harmonic space of a (pseudo-)hyperk\"ahler manifold} \label{sect3.2}
The {\it harmonic space} of a $4n$-dimensional  hk manifold $(M, g, J_\a)$ is   the trivial bundle   $\HM = M\times \SL_2(\bC) \to M$,  endowed with an integrable complex structure $\bI$ defined as follows. For each point $(x, U) \in \HM$ consider  the natural direct sum decomposition
$T_{(x, U)} \HM =  T_x M + T_U \SL_2(\bC) \simeq  T_x M +  \sl_2(\bC)$,  where  $T_U \SL_2(\bC)$ 
is identified with $\sl_2(\bC)$ by means of right invariant vector fields. 
Then let $\bI_{(x,U)}$ be the unique complex structure  on  $T_{(x, U)} \HM$ given by 
\beq
\bI_{(x,U)}|_{T_x M} \= I^{(z)}|_x,\  z = U{\cdot}[0:1] \in \bC P^1 (\simeq S^2) \ \text{and}\  \bI_{(x, U)}|_{\sl_2(\bC)} = \Jst\ ,
\label{sl2cxstr}\eeq
where  $\Jst$ is the  complex structure of $\sl_2(\bC)$ given  by the  multiplication by $\left(\smallmatrix i & 0\\ 0 & i\endsmallmatrix	\right)$. From the above  identification $T_U \SL_2(\bC) \simeq \sl_2(\bC)$,  it follows that along each fiber  $\{x\} \times \SL_2(\bC)$, the $\bI$-holomorphic distribution is generated by  right invariant   vector fields of $\SL_2(\bC)$.\par
 The collection  $\bI$ of such pointwise defined   complex structures is
 a globally defined almost complex structure on  $\HM$, which can be seen to be   integrable as follows.  
 The family of  restricted  complex structures  $\bI|_{M \times \{U\}} = I^{(z)}|_{M \times \{U\}}$ on the manifolds $ M \times \{U\}$, $U \in \SL_2(\bC)$,   is invariant under the natural left action of $\SL_2(\bC)$ on $\HM$. 
 Thus, the Lie derivative of an $\bI$-holomorphic vector field that is tangent to the (horizontal) slices $ M \times \{U\}$  by means of an infinitesimal transformations of this $\SL_2(\bC)$-action  always gives another  $\bI$-holomorphic vector field, which  is horizontal.  On the other hand, the infinitesimal transformations of the left action of $\SL_2(\bC)$ on each vertical fiber $\{x\} \times \SL_2(\bC) \simeq \SL_2(\bC)$ are nothing but the right invariant vector fields of $\SL_2(\bC)$. This implies that the Lie bracket  between a horizontal $\bI$-holomorphic vector field and a vertical $\bI$-holomorphic vector field is a horizontal $\bI$-holomorphic vector field. This property together with the fact that both the horizontal and vertical   $\bI$-holomorphic distributions are involutive proves that the whole $\bI$-holomorphic distribution of $\HM$ is involutive, i.e. that $\bI$ is globally integrable.\par
 \par
\smallskip
\begin{rem} \label{rem31} Note that $\HM$, considered as a principal $\SL_2(\bC)$-bundle over $M$,  can be (locally) identified  with a bundle of vertical complex  frames $(e_+, e_-)$ for the  fibers of the rank $2$ holomorphic vector bundle $\pi: {\bf H} \to M$ introduced  by Salamon in \cite{Sa} for general quaternionic K\"ahler manifolds. (Note that the fibers of ${\bf H}$ are the spaces ${\bf H}_x$ defined in \S \ref{adapted}.)  The explicit  construction of the complex structures of $\HM$ and ${\bf H}$   directly yields that the projection  
$$p: \HM \to {\bf H} \setminus \{0\}\ ,\qquad p((x, (e_+, e_-)) = e_+ \in {\bf H}_x,\  x \in M$$
is holomorphic. Thus $\HM$   fibers  holomorphically  over the twistor space $Z(M) = P({\mathbf H}) \simeq M \times \bC P^1$  with typical fiber given by  the subgroup  $B \subset \SL_2(\bC)$  of upper triangular matrices. Thus, $Z(M) \simeq \HM/B = M \times  \SL_2(\bC)/B$.  
\end{rem}
\par
By construction,  the harmonic space $\HM$ is  equipped with the $\bI$-invariant integrable distribution $\cD$ given by the tangent spaces of the leaves $M \times \{U\}$. 
The  two complex subdistributions of $\cD^\bC$, spanned by the  holomorphic and anti-holomorphic vector fields, will be denoted by 
$\cD_+$, $\cD_- \subset \cD^\bC$. 
By the remarks at the end of \S \ref{section321}, for  $(x, U) \in \HM$ there is at least one adapted frame 
$u: \bH^n \to T_x M$, more precisely,  a frame with  $u(\mathbf i) = I^{(z)}_x$,  $z = U{\cdot}[0:1]$, such that the corresponding  
 $2n$-tuples $(e_{+a})$ and $(e_{-a})$ are bases for the vector spaces  $\cD_+|_{(x, U)}$ and $\cD_-|_{(x, U)}$, respectively. 
 \par
\subsection{The complexified harmonic space} \label{sect3.3}
Consider an  $n$-dimensional complex manifold   $(N, J)$  and denote by  $\cA_J$ the complete atlas of holomorphic coordinates, i.e., of systems of coordinates  $\xi = (z^i): \cU \subset N \to \bC^n$ in which 
 the integrable complex structure $J$ has the standard form  $J = i \frac{\p}{\p z^j} \otimes d z^j - i  \frac{\p}{\p \overline z^k} \otimes d \overline z^k$. 
We define the {\it complexification} of $(N, J)$  as the pair $(N^{\bC}, \imath)$ given by:
 \begin{itemize}
\item[a)]  the complex manifold  $N^{\bC} \= N \times N$ having  complex  structure $\wt J$  defined at each point $(x, y) \in N\times N$ by $\wt J_{(x,y)} (v, w) {\=} J_x(v) {-} J_y(w)$
\item[b)] the standard  diagonal embedding 
$\imath: N \to N^{\bC}$,  $\imath(x) = (x, x)$.
\end{itemize}
Note that  the complex structure $\wt J$ of $N^\bC$ is  defined in such a way  that the corresponding  atlas $\cA_{\wt J}$  of  holomorphic coordinates is generated by coordinates of the form
$\wt \xi = (z^i, \overline{z'{}^j}): \cU \times \cV \to \bC^{2n}$ for  
some   (local) holomorphic coordinates  $(z^i)$,  $(z'{}^j)$ of $(N, J)$. 
So, the anti-holomorphic involution 
$\t: N^\bC \to N$,  $ \t(x,y) = (y,x)$ has a fixed point set which is precisely the totally real submanifold 
$\imath(N)  \simeq N$.
 \par
\smallskip
The above construction  yields  the following very convenient extension of  harmonic spaces. Let $\HCM$ 
be  the cartesian product  $\HCM = M \times M \times \SL_2(\bC)$,  equipped with the unique (integrable) complex structure $\bI^{(\bC)}$, which coincides with the right invariant  complex structure  along the  (vertical) leaves $\{x\} \times\{y\} \times \SL_2(\bC) \simeq \SL_2(\bC)$  (see  \eqref{sl2cxstr}) and with the complex structure  of the 
complexification  of $(M, I^{(z)})$, $z = U{\cdot}[0:1]$, along each horizontal leaf  $M \times M \times \{U\}$. In other words, $\HCM$ is  the union of the complexifications of the  manifolds $(M, I^{(z)})$, $z \in S^2$, and not the complexification of the harmonic space $\HM$.  Nevertheless, we call $\HCM$ the   {\it complexified harmonic space}. \par
\smallskip
Both of the  plus and minus distributions $\cD _\pm \subset T^\bC \HM$ naturally extend to  holomorphic  distributions on $\HCM$. In order  to see this, just  consider  the  distribution $\cD'$ on $\HCM$, determined by  the tangent spaces to  all leaves $M \times M \times \{U\}$  and the associated complex subdistribution  $\cD^{10}_{\bI^{(\bC)}} \subset \cD'{}^\bC$  spanned    by  $\bI^{(\bC)}$-holomorphic   vector fields. 
The  subdistributions  of  $\cD^{10}_{\bI^{(\bC)}}$ that project isomorphically onto the tangent spaces 
either on the first or the second copy of $M$ coincide with  the distributions $\cD_{\pm}$ at the points of the  real 
submanifold $\HM \subset \HCM$.  For simplicity, we shall denote these subdistributions of  $\cD^{10}_{\bI^{(\bC)}}$  also by $\cD_{\pm}$. 
\par
\section{Instantons on hk manifolds and  prepotentials}
\setcounter{equation}{0}
\subsection{Instantons on hyperk\"ahler manifolds} \label{definitioninstantons}
As usual,  let $(M, g, J_\a)$ be an hk manifold and denote by  $T^\bC_x  M \simeq \sH_x \otimes \sE_x,\ x \in M$,  the isomorphisms  described in \S \ref{adapted}.  The  space of complex $2$-forms $\L^2 T_x^{*\bC} M$ splits  into three irreducible $\SL_2(\bC){\cdot}\Sp_n(\bC)$ moduli:
 \beq \label{decomp} 
 \L^2 T_x^{*\bC} M  \simeq \bC \o_{\sH_x} \otimes S^2 \sE^*_x + S^2 \sH^*_x \otimes \bC\o_{\sE_x} + S^2 \sH^*_x \otimes \L^2_0 \sE^*_x\ .\eeq
Here $\o_{\sH_x}$ and $\o_{\sE_x}$ are the $\SL_2(\bC)$- and $\Sp_n(\bC)$-invariant symplectic forms of $\sH_x$ and $\sE_x$, respectively, and   $ \L^2_0 \sE_x$ is the irreducible $\Sp_n(\bC)$-submodule of $\L^2 \sE_x$  complementary to $\bC \o_{\mathbf E_x}$. Since the isomorphism  $T^\bC_x  M \simeq \sH_x \otimes \sE_x$
is unique up to an action of an  element in $\Sp_1{\cdot} \Sp_{p,q}$, the decomposition \eqref{decomp} is independent of the  isomorphism chosen.\par
Now,  given a gauge field  $(E, D)$ on an hk manifold, we split   the ($\bC$-linear extension of the)  curvature tensor 
$F_x,\  x \in M$   as follows:
 \begin{align*} 
 F_x = F^{(1)}_x + F^{(2)}_x\quad
  \text{with}\quad F_x^{(1)} &\in \bC \o_{\sH_x}{\otimes} S^2 \sE^*_x \otimes \End(E_x)\  \\
  F_x^{(2)} &\in (S^2 \sH^*_x {\otimes}(\bC \o_{\sE_x}+ \L^2_0 \sE^*_x))\otimes \End(E_x)\, . 
  \end{align*}
 A  gauge field  $(E, D)$ is called  {\it instanton} if the  $F^{(2)}$ component  of  its curvature tensor  vanishes
 everywhere.  Such instantons provide  minima  of
 the Yang-Mills functional $ \int_M |F|^2 \o_g$ and are thus solutions of the Yang-Mills equations \cite{ward84}.
Such instanton equations have been been studied by several authors \cite{Sa, MS, cdfn, Ti, ACD,De}. \par

Notice  that the vanishing   of $F^{(2)}$   corresponds  to simple conditions on the   components 
$F_{\pm a |\pm b} = F(e_{\pm  a},  e_{\pm b})$  with respect to   the  complex  frames $(e_{\pm a})$ defined in  \eqref{25}. In fact,  $F^{(2)} =  0$ if and only if 
\beq\label{3.10} F_{+a| + b} = F_{-a| -b} = 0\ ,\qquad F_{+a| - b} = - F_{-a |+ b}\ ,\qquad  F_{+a| - b} = F_{+b| -a}\ .\eeq
In four dimensions these are precisely the well-known anti-self-duality equations.
\par
\subsection{Central and exponential-central gauges  on harmonic spaces}\label{centr-exp}
 A  gauge   $\f = (\f^\cV, \f^G)\colon P|_{\cV} {\to} \cV \times G$  for  the $G$-bundle $P$ naturally corresponds to a gauge  for its lift $P'$ on the harmonic space $\HM$,  namely the gauge $\f$ defined on  the restriction of $P'$ to $\cV \times \SL_2(\bC) \subset  \HM$  by 
 $$\f : P'|_{\cV \times \SL_2(\bC)}  \to \cV \times  \SL_2(\bC) \times G\ ,\quad \f(u, U) \= (\f^\cV(u), U, \f^G(u))\ .$$
 Such a gauge  is called the {\it central gauge  determined by  $\f$} \cite{gios_book, ACD}. \par
\smallskip
Let us now  introduce a very convenient special class of  central gauges.
Given $x_o \in M$, for any unit vector $v \in T_{x_o} M$ we denote by  $\g^{(v)}_t = \exp_{x_o}(t v)$
the radial geodesic determined by $v$. Now let $\cV \subset T_{x_o} M$ be  a neighbourhood of $x_o$
such that  the exponential  map $\exp_{x_o}: \cV \subset T_{x_o} M \to M$ is a diffeomorphism onto its
image $\cU = \exp_{x_o}(\cV)$.
A gauge $\f: P|_{\cU} \to \cU \times G$ on the domain  $\cU = \exp(\cV)$ is called {\it exponential} if the corresponding potential $A$ satisfies the following conditions:
\begin{itemize}
\item[(a)] $A|_{x_o} = 0$ and 
\item[(b)] $A(\dot \g^{(v)}_t) = 0$ for  all    vectors $\dot \g^{(v)}_t$  tangent  to  the radial geodesics.
\end{itemize}
We shall call  the central gauges for $P'$ on $\HM$ determined by    exponential gauges for $P$ on $M$
{\it exponential central} (or just {\it exp-central}).
 We recall  that  for  any  $x_o \in M$  there is always   an exponential gauge on some appropriate neighbourhood of  $x_o$ \cite{Uh1}. Thus, {\it  for any $(x_o, U) \in \HM$ there exists an exp-central gauge for $P'$ on a neighbourhood of $(x_o, U)$}.\par
It is also known that  for  all cases in which  $(E, D)$ is the complexification of a gauge field $(E^o, D^o)$ with compact structure group   $G^o \subset G$,   if  $A$ is the potential  in an exponential gauge $\f$ for the bundle $P^o$  with domain $\cU\subset M$, then there exists  a constant $c_\cU$, which   depends only on $\cU$, such that 
\beq\label{Uhlenbeck-central} \|A \|_{\cC^0(K, \gg)} \leq c_{\cU} \|F\|^{(\f)}_{\cC^0(K, \gg)}\quad \text{for} \ K \subset\subset \cU\ ,\eeq
 (see e.g. \cite[Lemma 2.1]{Uh1}). Similar  estimates clearly   hold for  potentials and curvatures of the lifted gauge fields $(E', D')$  in exp-central  gauges.\par
\par
\subsection{Prepotentials for instantons on hk manifolds}
\label{sectprepotentials}
 We recall that {\it any hk manifold  $(M, g, (J_\a))$  has a natural structure of a  real analytic manifold and that in such a structure,    the tensors $g$  and $J_\a$ are  real-analytic} \cite{Le}. Hence, if we  lift  $g$ to $\HM$ as a tensor field with values in $\cD^* \otimes \cD^*$, using  real-analyticity,  for each point $(x, U) \in \HM \subset \HCM$  we may determine a tubular $\SL_2(\bC)$-invariant neighbourhood $\cW \subset \HCM$ of $(x, U)$, to which  $g$ extends as a $\bC$-linear tensor field  in  $\cD^{\bC*}\times \cD^{\bC*}$. \par
 We claim that an analogous  extension property holds also for an instanton on an  hk manifold  provided that it    is a complexification of an  instanton with compact structure group. To  prove this,  let us introduce some additional convenient notation.  Given a gauge field $(E, D)$ associated with $(P, \o)$,  let   $\gH \subset TP'$ and $\gH^\bC \subset T^\bC P'$ be  respectively the {\it real} and  {\it complex} horizontal distributions of $P'$  given by the kernels of the lifted connection  $\o'$ on $P'$.  
 Further, for  any (real or complex) vector field $X$ on $\HM$, let us   denote by  $X^h$ the uniquely associated   vector field on $P'$ with values  in  $\gH$ or $\gH^\bC$ which  projects  onto $X$. 
 
 \begin{prop} \label{prop41} If $(E, D)$ is an instanton on $(M, g, (J_\a))$, which  is the complexification of a gauge field with compact structure group $G^o$,   all data of the lifted pairs $(E', D')$ and   $(P', \o')$  on $\HM$ are  real analytic. Moreover,  there is a complex structure $J'$ on $P'$,  invariant under $G = (G^o)^\bC$, which  
makes $p': P' \to \HM$ and the associated vector bundle $p^{E'}: E' \to \HM$  holomorphic bundles over $\HM$   with  $J'$-invariant connections. \par
 It follows that, for any $(x_o, U) \in \HM \subset \HCM$, there are  unique   real-analytic extensions  of  $(E', D')$ and   $(P', \o')$   to some    $\SL_2(\bC)$-invariant tubular neighbourhood $\cW \subset \HCM$ of the real submanifold $\cU \= \cW \cap \HM$ containing   $(x_o, U)$.
 Further, both the extended bundles $E'$, $P'$ have naturally extended complex structures  that  make them   holomorphic bundles over $(\HCM, \bI^{(\bC)})$  and the extended connection $\o'$ is   $J'{}^\bC$-invariant with respect to the complex structure $J'{}^\bC$ of $P'$. \par
 Finally, the components   in holomorphic gauges of the extended connection $\o'$  and of its $(1, 0)$-potential are holomorphic functions of any set of  complex coordinates $(z^\ell, w^m, (u^i_r))$ of the complex manifold $(\HCM = M \times M \times \SL_2(\bC), \bI^\bC)$ that correspond to complex coordinates $(z^\ell, \overline z^m = w^m, (u^i_r))$ of the submanifold    $\HM \subset \HCM$, 
 \end{prop}
\begin{pf}
To prove these statements,  we   need the following simple  lemma.

\begin{lem} \label{char} The complex gauge field $(E, D)$   is an  instanton if and only if 
the curvature $F'$ of its lift $(E', D')$ on $\HM$ is such that 
\beq \label{charact}   F'(X_+, Y_+)= 0 =    F'(\overline{X_+}, \overline{Y_+}) \qquad \text{for}\ \ X_+, Y_+ \in \cD_+\ . \eeq
\end{lem}
\begin{pflemma}  As observed above,  each  space $\cD_+|_{(x,U)} $ or (more precisely, its isomorphic  projection  onto $T^\bC_x M$) is spanned by the  vectors   $ e_{+a}$  of the adapted frames $u: \bH^n \to T_xM$. Hence  the curvature  of $(E,D)$ satisfies the first  condition in \eqref{3.10} 
if and only if  $F'(X_+, Y_+)= 0$  and $F'(X_-,  Y_-)= 0$ for  $X_\pm $, $Y_\pm \in \cD_\pm$. The other  conditions in \eqref{3.10}  are  direct consequences of the decomposition of $\cD^\bC|_{(x,U)} \simeq T^\bC_x M$  into irreducible $\SU_2$-moduli.  \end{pflemma}

Due to this lemma and  \eqref{2.4},   the curvature  $2$-form $\O'$ of the lifted connection $\o'$ identically vanishes on any  pair of  horizontal lifts $X^{01h}, Y^{01h} \in \gH^\bC$ of vector fields $X^{01},  Y^{01}$ in the anti-holomorphic distribution $\cD_- + V^{01} = T^{01} \HM$ of $(\HM, \bI)$, where $\cD_-$ and $V^{01}$ are the $\bI$-anti-holomorphic horizontal and vertical distributions  of $\HM$ described  in \S \ref{sect3.2}.  It follows that  the subbundle $\cS^{01} \subset T^\bC P'$ generated by the vectors  $X^{01h}$ and the   anti-holomorphic  vertical  distribution of $P'$  is an  involutive complex distribution. The same holds for the  subbundle $\cS^{10}  = \overline{\cS^{01}}\subset T^\bC P'$.  Hence  the direct sum decomposition 
$T^\bC P' = \cS^{10} + \cS^{01}$   corresponds to  a $G$-invariant  integrable complex structure $J'$ on $P'$. Consequently,  there is an atlas of  complex charts  for the complex manifold $(P', J')$ making $P'$ a holomorphic $G$-bundle over $\HM$. Moreover, 
the lift $P^{o}{}'$ on $\HM$ of the bundle $P^o$  with compact structure group is necessarily a real analytic submanifold of $P'$ since it is the  fixed point set of  an appropriate real analytic involution. One can also check that the restricted distribution  $\gH|_{P^o{}'}$ coincides with the distribution given by the  $J'$-invariant subspaces of the tangent spaces of $P^{o}{}'$. Since the latter is a real analytic distribution on $P^o{}'$ and   the  distribution $\gH$ of $P'$ is the unique $G$-invariant extension of $\gH|_{P^o{}'}$, we conclude that also $\gH$ is $J'$-invariant  on $P'$.  Consequently, the first claim follows immediately.\par
Concerning the second claim, the same arguments as above yield   the existence of a  complex structure  ${J'}^\bC$ on the extended bundle $P'$ (and, consequently, a corresponding complex structure on the associated bundle  $E'$), which  makes it a holomorphic bundle over the complex manifold  $(\HCM,  \bI^{(\bC)})$ and   leaves  invariant the horizontal distribution $\gH$ determined by the extended connection $1$-form $\o'$.    
The final claim is a  consequence  of  the fact that the extension to $\cW \subset \HCM$ of every real analytic datum on $\HM$ 
is obtained by considering the power series of such a datum in the variables $(z^\ell, w^m \= \bar z^m, (u^i_\ell))$ as a power series in the independent complex variables $(z^\ell, w^m, (u^i_\ell))$, thus holomorphic in both the variables $z^\ell$ and $w^m$. The holomorphy in the   $u^i_\ell$  follows from the fact that  $\o'$ is  the lift  of a  connection form of  $p: P \to M$ to the bundle  $p': P' \to \HM$. 
\end{pf}
Consider   the  basis  $(H^o_0, H^o_{++}, H^o_{--})$ of $\sl_2(\bC)$ defined by 
\beq\label{Hzero}  H^o_0 \= \begin{pmatrix} 1 & 0\\ 0 & - 1 \end{pmatrix}\ ,\quad H^o_{++} \= \begin{pmatrix} 0 & 1\\ 0 & 0 \end{pmatrix}\ ,\quad H^o_{--} \= \begin{pmatrix} 0 & 0\\ 1 & 0 \end{pmatrix}\eeq
and denote by  $H_0$, $H_{++}$, $H_{--}$ the  associated   holomorphic vector fields on $\HM$, determined as  holomorphic parts of   the infinitesimal transformations  of the right actions of  the one-parameter groups generated by  $H^o_0, H^o_{++}, H^o_{--}$.  The restrictions of such vector fields $H_\a$ to each vertical fiber $\{x\} \times \SL_2(\bC) \simeq \SL_2(\bC)$ are   left invariant and generate a Lie algebra isomorphic to $\sl_2(\bC)$.\par
We may now   prove the theorem, on which all our results  are based (see also \cite{ACD},  Prop. 7 and Thm.4). 
\begin{theo} \label{thm53}  Given  a real analytic complex  instanton  $(E, D)$  and  a point $(x_o, x_o, U) \in \HM$, there exists an $\SL_2(\bC)$-invariant neighbourhood $\cU \subset \HCM$ of $(x_o, U)$ and a  holomorphic gauge 
$\varphi: P'|_{\cU} \to \cU \times G$ of the $\HCM$-extension  of the lifted bundle  of $P$ with   associated $(1,0)$-potential $A$ satisfying the conditions
\beq \label{claim1} A(H_0) =  A(X_-)  = 0\ ,\ X_- \in \cD_-\ .
\eeq
Such a potential $A$ is uniquely   determined by the $\gg$-valued function  $A_{--} \= A(H_{--})$ in the following sense: given the lift  $(\wh E', \wh D')$  on $\HM$ of a  real analytic complex  instanton $(\wh E, \wh D)$ on $M$,  if    $\wh E'|_{\cU}$ coincides with  $E'|_{\cU}$ and  if furthermore the  $(1,0)$-potential  $\wh A$ of   $(\wh E', \wh D')$ in a holomorphic gauge  satisfies
\eqref{claim1},  then 
\beq \wh D'|_{\cU \cap \HM} {=} D'|_{\cU \cap \HM}  \ \Longleftrightarrow \  
 \label{claim2} \wh A(H_{--})|_{\cU \cap \HM} {=} A(H_{--})|_{\cU \cap \HM}\ .\eeq
\end{theo}
\begin{rem} Note that  by $ \wh D'|_{\cU \cap \HM} {=} D'|_{\cU \cap \HM} $ we  mean that the  differential operators $D'$, $\wh D'$  are {\it identical},  not just equivalent up to an automorphism (gauge transformation) 
of   $\wh E'|_{\cU}$. 
However,  \eqref{claim2} implies a  bijective  correspondence   between covariant derivatives and their  potential components $A_{--} = A(H_{--})$. This  induces  a natural  bijective correspondence between  equivalence 
classes of  instantons   up to  local automorphisms  and equivalence  classes  of    $(1,0)$-potential components  $A_{--} = A(H_{--})$ up to  appropriate gauge transformations  \eqref{changepot}.  
In fact, the required gauge transformations   are precisely  those leaving condition  \eqref{claim1} unchanged. 
\end{rem}
\begin{pf}  Let us extend  $(P', \o')$ to an   $\SL_2(\bC)$-invariant tubular neighbourhood $\cW  \subset \HCM$  of $\HM$ and consider the  two $G$-invariant   complex distributions  $\cD^h_-$ and  $\cD^h_- \oplus < H^h_0> $ on the  extended bundle $P'$, 
   the former generated   by the horizontal lifts $X^h$ of the complex  vector fields  in $\cD_-$, the latter  generated by the complex vector fields  in $\cD^h_-$ and  the horizontal lift $H^h_0$ of the  holomorphic  vector field  $H_0$ of $\HCM$.  We recall that  the distribution $\cD_-$ is spanned by the complex vector fields $e_{-a}$ described in \S \ref{adapted}.  Using  this and the   identification, described in Remark \ref{rem31}, between the points of  $\HM$ and the frames for the fibers ${\bf H}_x$ of the bundle ${\bf H}$,  the distribution $\cD_-$ is seen
   to be invariant under the flow of  the vector field $H_0$. This implies that  the horizontal distribution $\cD^h_-$ is invariant under the flow of the horizontal vector field  $H^h_0$. \par
By  Lemma \ref{char}, the curvature  $2$-form $\O'$ of $\o'$  vanishes identically on  the distribution $\cD^h_-$, meaning that    $\cD^h_-$  is involutive.  Moreover, from the last claim in Proposition \ref{prop41}, $\cD^h_-$ is generated by holomorphic vector fields of the holomorphic bundle $\pi: P' \to \cW \subset \HCM$.
Hence, by the complex Frobenius Theorem, for each $y_o  = (x_o, x_o, U, g)\in  P'|_\HM$  there is an 
$\SL_2(\bC) {\times} G$-invariant neighbourhood $P'|_{\cU}$ of $y_o$, which is holomorphically foliated by integral leaves of $\cD^h_-$. Note that the  union of  the $H^h_0$-orbits of  the points of one  such integral leaf is  an integral leaf of the larger complex distribution $\cD^h_{-}   \oplus < H^h_0> $. It follows   that  $P'|_{\cU}$ is actually holomorphically foliated by the integral leaves of this  larger distribution and   we may  consider   a holomorphically parameterised  family of   integral leaves of  this distribution, which fills a complex submanifold $\cS'$ 
 transversal to the $G$-orbits at each point.   
Without loss of generality, we may also assume that  such a submanifold  $\cS'$ projects biholomorphically  onto  an $\SL_2(\bC)$-invariant  neighbourhood $\cU \subset \HCM$ of $(x_o, x_o, U)$ and hence it is a graph of a holomorphic section  of the $G$-bundle $P'$.  Associated with  such  a section, there is a  unique holomorphic gauge   $\f: P'|_{\cU} \to \cU \times G$  mapping $\cS'$  onto the trivial section   $\cU \times \{e\} $ of  the trivial bundle $ \cU \times G$.  By construction,  the  $(1,0)$-potential $A$  in such a  gauge   satisfies   \eqref{claim1}.
\par
\smallskip
Suppose now that  $(\wh E', \wh D')$  is a real analytic  lifted  instanton  associated with the pair  $(\wh P', \wh \o')$,  and assume that $\wh A(H_{--})|_{\cU \cap \HM} {=} A(H_{--})|_{\cU \cap \HM}$. By real analyticity,   we may   assume  that $\cU$ is such that   $\wh A(H_{--})|_{\cU} = A(H_{--})|_{\cU}$. 
Thus,  the horizontal lifts $H^h_\a$, $\wh H^h_\a$ of the  holomorphic vector fields $H_\a$, $\a \in \{0, ++, --\}$, determined by the connection forms  $\o'$ and $\wh \o'$, respectively,  have the forms
\begin{align}
\nonumber &\wh H^h_0 = H^h_0\ ,\qquad \wh H^h_{--}  = H_{--}  + \wh A_{--} =  H_{--}  +  A_{--} =  H^h_{--}   ,\\
\label{5.14} &\wh H^h_{++}  = H_{++}  + \wh A_{++}\ ,\qquad  H^h_{++}  = H_{++}  + A_{++}\ ,
\end{align}
where we denote by $A_{\pm\pm} \= A(H_{\pm\pm})$,  $\wh A_{\pm\pm} \= \wh A(H_{\pm\pm})$  the components along the vector fields $H_{\pm\pm}
$ of the corresponding   $(1,0)$-potentials.   Setting $B_{++}  \= \wh A_{++}  - A_{++}$, 
the expansions \eqref{5.14}  can be written as  
\beq \label{rel1}  \wh H^h_0 = H^h_0\ ,\qquad \wh H^h_{++} = H^h_{++} + B_{++}\ ,\qquad \wh H^h_{--} = H^h_{--}\ .\eeq
Since the considered principal bundles $P'$,  $\wh P'$ are lifts to $\HM$ of bundles over $M$,  the Lie brackets  among the horizontal lifts of  vector fields $H_A$ coincide with  the Lie brackets among their projections (see \eqref{2.8})  so that 
$$[\wh H^h_0, \wh H^h_{++}] =  2 \wh H^h_{++}\ ,\qquad [\wh H^h_{++}, \wh H^h_{--}] = \wh H^h_0\ . $$
This and \eqref{rel1}  imply that the  $\gg$-valued holomorphic function  $B_{++}  $ is such that 
$$H^h_0{\cdot} B_{++} = 2 B_{++} \ ,\qquad H^h_{--}  {\cdot} B_{++}  = 0\ .$$ 
On the other hand, by Proposition \ref{prop28},   the vector fields $H^h_0$, $H^h_{++}$,  $H^h_{--}$ generate a right holomorphic  $\SL_2(\bC)$-action on $P'|_{\cU}$ and  the bundle $P'|_{\cU}$ is foliated by regular orbits of this action  (thus, each such orbit is identifiable with a copy of $\SL_2(\bC)$).   We may therefore apply  Lemma 5.3 in \cite{DS}  (see also \cite{gios_book}, \S 4.3) -- considered for  the pair of vector fields $H_0$,  $H_{--}$ in place of the  pair $H_0$, $H_{++}$ -- to the restrictions of $B_{++}$ to each such orbit.  
This lemma implies  that  $B_{++}$ is  identically vanishing along each such orbit. It follows  that  $B_{++} \equiv 0$, i.e. $\wh A_{++} \equiv A_{++}$, and   that $\wh H^h_{++}  = H^h_{++}$. \par
\smallskip
In order to conclude  that $\wh D' = D'$,  it is now sufficient to prove  the existence  of  a collection of holomorphic  vector fields $(e_{+a},  e_{-b})$ on $\cU \subset \HCM$ with the following two properties: a) they span a distribution  which is  complementary to the one generated   by the vector fields  $H_A$ and b)  the remaining   components  $\wh A_{\pm a} \= A(e_{\pm a})$ and $A_{\pm a} =  A(e_{\pm a}) $ of the potentials of  $\wh D'$ and $D'$ are identical.  Let us consider   a  set of (locally defined) holomorphic vector fields  $(e_{-a})$ generating  the distribution $\cD_-$ and  projecting   pointwise onto  vectors $e_{-a} \in T^\bC M$  determined by adapted frames    as  in \S \ref{adapted}. Since the  vector fields $e_{-a}$ are eigenvectors with eigenvalue $-1$ for the action of the complex vector field $H_0$, the vectors 
$e_{+a} \= [H_{++}, e_{-a}]$ are eigenvectors with eigenvalue $+1$ for the same action. This implies that the   $e_{+a}$ generate the distribution $\cD_+$ and that    $(e_{+a}, e_{-b})$ is a collection of generators for  the distribution $\cD^{10}_{\bI^\bC}$ complementary to the distribution spanned by the $H_A$. Moreover, their corresponding lifts  $e^h_{\pm a}$, $\wh e_{\pm a}^h$, determined by  the two connection forms $\o'$,  $\wh \o'$,  are such that  
\beq \label{rel2} [H^h_{++},  e^h_{-a}] = e^h_{+a}\ ,\qquad[\wh H^h_{++},  \wh e^h_{-a}] = \wh e^h_{+a}\ .\eeq
 Setting $A_{\pm a} {\=} A(e_{\pm a})$, $\wh A_{\pm a} {\=} \wh A(e_{\pm a})$ and  recalling that, by hypothesis,  
 $A_{-a} = \wh A_{-a} = 0$, 
  we have that  \eqref{rel2}  implies  $ - e^h_{-a}{\cdot} A_{++} = A_{+a}$ and  
   $- e^h_{-a}{\cdot} \wh A_{++} = A_{+a}$.
Since we have proven that $\wh A_{++}   \equiv  A_{++}$,   this  gives  $ \wh A_{\pm a}   \equiv  A_{\pm a}$, as  desired. 
\end{pf}
We call the   $\gg$-valued map $A_{--}|_{\cU \cap \HM} \=  A(H_{--}|_{\cU \cap \HM})$,
which uniquely determines the extended  holomorphic function  $A_{--}|_\cU$ and thereby the gauge field 
$(E'|_\cU, D')$,  a {\it  prepotential on $\cU \cap \HM$ for  the  instanton $(E, D)$}.  The holomorphic gauges of the (extended) lifted bundle $P'$ in which the potential of $\o'$ satisfies  \eqref{claim1} are called {\it analytic}. \par
 \begin{rem} 
The previous literature on the harmonic space formulation e.g. \cite{GIO, gios, gios_book, ACD} used gauge conditions
$A_0 = A_{+a} = 0$ and prepotentials  $A_{++}$. Here we choose to reverse the role of the signs. This has the advantage
that  prepotentials  are holomorphic rather than  anti-holomorphic  with respect to the  complex structure $\bI$ of $\HM$
 (see Remark \ref{remark47}).
 \end{rem}
 \subsection{Analytic gauges,  bridges and normalisations}
 \label{normalisation}
Assume that our instanton $(E, D)$ is  the complexification of an  instanton with compact structure group $G^o$ over the hk manifold $(M, g, (J_\a))$. Around  each     $(x_o, U) \in \HM$  there are two very important  classes   of  gauges to be  considered:  the  exp-central gauges and the analytic gauges. Let us briefly   compare their main features: 
\begin{itemize} [itemsep=2pt, leftmargin=18pt]
\item[--] If $\f: P'|_{\cV \times \SL_2(\bC)} \to \cV \times \SL_2(\bC) \times G$ is (the restriction to $\HM \subset \HCM$ of) an {\it analytic} gauge, the corresponding  holomorphic potential  $A$ for $\o$   is such that $A_0 \=A(H_0)$ is identically vanishing  and $A(X_{-}) = 0$  for any vector field $X_- \in \cD_-$. In contrast with this,    the functions $A_{\pm\pm} \= A(H_{\pm\pm})$ and  the functions $A(X_+)$,  $X_+ \in \cD_+$, are in general non-trivial.
\item[--] If $\wt  \f: P'|_{\cV \times \SL_2(\bC)} \to \cV \times \SL_2(\bC) \times G$ is an {\it exp-central} gauge,  
the  corresponding potential  $\wt  A$   is such that the components $\wt A_\a  = \wt A(H_\a)$,
 $\a \in\{0, ++, --\}$, identically vanish,  while  all  functions $\wt  A(X_\pm)$, $X_{\pm} \in \cD_{\pm}$, are in general non-trivial, being nevertheless constrained by  the  conditions  (a),  (b) in \S \ref{centr-exp}. 
\end{itemize}
The  maps $g: \cU \times \SL_2(\bC) \to G$,  which give the    gauge transformations $h \to g_{(x, U)} {\cdot} h$ from  central gauges  to  analytic gauges  are usually  called {\it bridges} (e.g. \cite{GIO, gios_book, ACD}). 
Now,   since  $G$ is the complexification of the compact Lie group $G^o$,  it  is reductive and the exponential  $e^{(\cdot)}: \gg \to G$ is a surjective local diffeomorphism. This means that any bridge $g_{(x, U)}$ can be  
written as  $g_{(x,U)} = e^{\psi_{(x, U)}}$ for some appropriate $\gg$-valued function $\psi$.  
We call  $\psi$  a {\it  $\gg$-bridge}.\par
\smallskip
In the next lemma we shall give the proof of existence of  bridges  and $\gg$-bridges, having  the special property that  the prepotentials determined in the newly built analytic gauges satisfy  additional {\it normalisation conditions},  which drastically reduce their degrees of freedom. Such normalised analytic gauges can be considered as  complex analogues of  the   Coulomb gauges for arbitrary gauge fields. \par
\smallskip
  In order to properly state such a normalisation, we  need to introduce some appropriate notation.  
  Given   $x_o \in M$ and a (sufficiently small) simply connected neighbourhood $\cV \subset M$ of $x_o$, 
  let    $(e_{+a}, e_{-b})$ be  a $4n$-tuple of holomorphic vector fields of $\HCM$ that generate  the distributions $\cD_{\pm} \subset T^\bC (\cV\times \cV \times \SL_2(\bC))$ ($ \subset T^\bC \HCM)$,   constructed as  in the proof of Theorem \ref{thm53}. Then, pick an element $U_o \in \SL_2(\bC)$,  
  say $U_o =\left( \smallmatrix  1  &  0\\  0  & 1\endsmallmatrix\right)$, 
and let  $\l^a, \mu: \cV\times \cV \times \SL_2(\bC) \to \bC^{2n}$, $a = 1, \ldots, 2n$, be  a set of  $2n+1$ holomorphic functions satisfying the following conditions
  \beq \label{problemino} 
\begin{split}  & H_{++}{\cdot} \l^a  = 0\ , \qquad  H_{++} {\cdot} \mu = -1\ ,\quad H_0{\cdot} \l^a =   
H_0{\cdot} \mu = 0\ ,\\
 & e_{+b}{\cdot}\l^a =  \d^a_b\ , \qquad e_{+b}{\cdot} \mu = 0\ , \qquad \l^a|_{(x_o, x_o, U_o)} =  \mu |_{(x_o, x_o, U_o)} = 0\ .
 \end{split}\eeq
By \cite[Lemma 5.2]{DS},  functions satisfying the first  line of these conditions   surely exist and are determined up to addition of a holomorphic  function constant along each $\SL_2(\bC)$-orbits.  Using the commutation relations between  $H_{++}$, $H_0$ and $e_{+a}$, we can see that also  the second line of  these conditions can  be  satisfied,  fixing   the $\l^a$ and $\mu$ completely. Now consider the   holomorphic distribution $ \wt \cD \subset  T^\bC (\cV\times \cV \times \SL_2(\bC))$ generated  by the vector fields 
$$e_{+a}\ ,\qquad H_{++} \  ,\qquad \wt H_{--} \= H_{--} + \l^a e_{-a} + \mu H_0\ , $$
  which  can be directly checked to be  involutive.
 Finally let $\wh S$  and $S$ be  two complex submanifolds of  $\cV \times \cV \times \SL_2(\bC) \subset \HCM$,   one included in the other, which are integral leaves  of the holomorphic distributions $\cD_+$ and $\wt \cD$, respectively, and  both  passing through $(x_o, x_o, U_o)$. Note 
   that each tangent space of $S$  is complementary  to  $\cD_-|_{(x, y, U)} + <H_0|_{(x, y, U)}>$. \par
   We may now state and prove the advertised existence result.
 \begin{lem}\label{lemmabridge}  Let $\wt \f: P'|_{\cV \times \SL_2(\bC)} \to \cV \times \SL_2(\bC) \times G$ be an exp-central gauge, determined by an exponential gauge for $(P, \o)$ on a simply connected  neighbourhood $\cV \subset M$ around   $x_o \in M$,  and  $\wt A$ the associated potential for $\o$.  Let also  $\wh S \subset S  \subset  \cV\times \cV\times \SL_2(\bC)$ be the   pair of complex submanifolds  passing through $(x_o, U_o)$ described above. 
If $\cV$ is sufficiently small,  there exists  an analytic gauge $\f: P'|_{\cV\times \cV \times \SL_2(\bC)} \to \cV \times \cV \times \SL_2(\bC) \times G$, 
    in which the  prepotential $A_{--}$  is such that 
 \begin{multline} \label{4.13} 
  H_0{\cdot} A_{--}   = - 2 A_{--}\ ,   \quad e_{-a}{\cdot} A_{--} = 0\ ,\quad H_{--}{\cdot} A_{--}  = 0 \ ,\\
   A_{--}|_{\wh S} = -  \l^a \wt A_{-a}|_{\wh S}\ , \quad H_{++}{\cdot} A_{--}|_{\wh S}  = -  \l^a \wt A_{+a}|_{\wh S} .
\end{multline}
\end{lem}
\begin{pf}  For simplicity, we  use  $\wt \f$ to  identify  $P'|_{\cV \times \SL_2(\bC)}$ with  $\cV \times \SL_2(\bC) \times G$ so that we may assume that the  considered   exp-central gauge  
is nothing but   the  identity  map. Let us consider  the integral leaves  in $\cV \times \cV \times  \SL_2(\bC)\times G$ of the holomorphic  distribution  $\cD^h_- + <H^h_0>$, which passes through
the points of  the manifold $S\times \{e\} \subset \cV \times \cV \times \SL_2(\bC) \times \{e\}$. Being horizontal, they are  transversal to the $G$-orbits and,  by dimension counting, they fill a submanifold  $\cS' \subset P'$  which  projects diffeomorphically onto an $\SL_2(\bC)$-invariant  neighbourhood  $\cV' \times  \cV' \times \SL_2(\bC)$ of $(x_o, x_o, U_o)$. Thus  there is a gauge $\f^o: P'|_{\cV \times \cV \times \SL_2(\bC)} \to \cV \times \cV \times \SL_2(\bC) \times G$, which maps $\cS'$ into the submanifold  $\cV\times \cV \times \SL_2(\bC) \times \{e\}$  and which  we may assume to satisfy the condition  $\f^o|_{S \times \{e\}} = \Id_{S \times \{e\}}$ (see also  the proof of Theorem \ref{thm53}).  By construction, $\f^o$ is an analytic gauge and  the bridge $g_{(x, y, U)}$ from  
$\wt \f = \Id$ to $\f^o$ is such that  $h \mapsto g^o_{(x, y, U)}{\cdot} h = e{\cdot} h$ 
for each $(x, y, U) \in S $. Hence, writing this  bridge in the form $g^o_{(x, y, U)} = e^{\psi^o(x, y, U)}$ for an appropriate $\gg$-bridge $\psi^o$, we have that  $\psi^o|_{S} \equiv 0$. Note  that  $\psi^o$ satisfies  the   equations $H_0{\cdot} \psi^o  = 0$ and $e_{-a}{\cdot} \psi^o  +  \wt A(e_{-a}) = 0$ because of  the following three properties: a)     in any central gauge the potential $\wt A$ satisfies  $\wt A(H_0) = 0$; , b) 
 the  potential $A^o$ in the analytic  gauge $\f^o$ satisfies   \eqref{claim1} and  c) the  potentials $\wt A$ and $A^o$ are related by   \eqref{changepot}.\par
 \medskip
 Let us now denote by $A^o_{--}$ the prepotential of the given instanton  in this gauge $\f^o$. Expanding the identities $\cF^{\f_o}(H_0, H_{--}) = \cF^{\f_o}(e_{-a}, H_{--}) = 0$ in terms of the potential $A^o$ in this gauge and recalling that, being an analytic gauge, we have $A^o_0 = A^o_{-a} = 0$, we find  that $H_0{\cdot}A^o_{--} = - 2 A^o_{--}$ and $e_{-a}{\cdot} A^o_{--} = 0$, i.e. the first two conditions of \eqref{4.13}. We remark that 
 {\it these two  conditions are  satisfied by  any prepotential}, being merely consequences of the above properties of the curvature. \par 
 \smallskip
 However, the  prepotential   $A^o_{--}$ does not necessarily satisfy   the other conditions in \eqref{4.13} also. 
 To get a prepotential with such additional properties,  we   need to further change  $\f^o$ into  a new  (further restricted) gauge $\f$, which  preserves the property of being analytic, i.e.  with $(1,0)$-potential components $A_0$ and $A_{-a}$  identically vanishing  and   with  $A_{--}$  satisfying the first pair of  equations in \eqref{4.13}.  In order to determine  this new analytic gauge $\f$, let us consider the $(2n+2)$-dimensional involutive distribution $\wh \cD \subset T^\bC (\cV \times \cV \times \SL_2(\bC) \times G)$  generated by the holomorphic vector fields $H_0$,  $H_{\pm\pm} + A^o_{\pm\pm}$  and $e_{-a}$. 
 Then, for each point $(x, y, U_o, g)$ of the  manifold $\wh S \times G$, consider the unique  integral leaf  
 through $(x, y, U_o, g)$ of this distribution 
 $\cT^{(x,y, U_o, g)} \subset \cV \times \cV \times \SL_2(\bC) \times G$. Now we may  determine a $G$-equivariant $\gg$-valued holomorphic function $h_{--}$ on $\cT^{(x,y, U_o, g)} $  (the $G$-equivariance being with respect to the standard  right action of $G$ on $\cV \times \cV \times \SL_2(\bC) \times G$ and the adjoint $G$-action on $\gg$) satisfying the equations
\beq \label{original} H_0{\cdot}  h_{--}= - 2 h_{--} \ ,\ \  H_{--}{\cdot}  h_{--} = - H_{--}{\cdot} A^o_{--} -  [h_{--}, A^o_{--}]\  , \ \  e_{-a}{\cdot} h_{--} = 0\ .
\eeq
Once again, the existence of such an $h$ can be checked using  \cite[Lemma 5.3]{DS}. Indeed,  one can construct a solution to the first two conditions as follows. Along the image of   $\cT^{(x, y, U_o, g)}$ under  the inverse gauge transformation $\wt \f \circ \f^{o-1} (= \f^{o-1})$  (note that, by $G$-equivariance of the trivialisation,  each intersection of the   submanifold
$ \f^{o-1}(\cT^{(x,y, U_o, g)} ) $ with a vertical set $\{(\wt x,\wt y)\} \times \SL_2(\bC) \times G$  is entirely included in a submanifold of the form  $\{(\wt x, \wt y)\}  \times \SL_2(\bC) \times \{g\}$) 
 we may consider a $G$-equivariant $\gg$-valued solution $\wt h_{--}$  to  the differential problem 
\beq \label{original1} H_0{\cdot}  \wt h_{--}= - 2 \wt h_{--} \ ,\quad H_{--}{\cdot}  \wt h_{--} = - H_{--}{\cdot} \left(A^o_{--}\circ \f^o \right) \ , \quad e^o_{-a}{\cdot} \wt h_{--} = 0\ .\eeq
By \cite[Lemma 5.3]{DS},  applied to the pair of vector fields $H_0$,  $H_{--}$ in place of the  pair $H_0$, $H_{++}$, such a solution $\wt h_{--}$ exists. The corresponding  function $h_{--} = \wt h_{--} \circ \f^{o-1}$  is therefore  a  solution to \eqref{original}. Note that the  solution $\wt h_{--}$  to \eqref{original1} is uniquely  determined up to addition of a solution $\wt k_{--}$  of the associated  system 
$$ H_0{\cdot}  \wt k_{--}= - 2 \wt k_{--} \ ,\qquad H_{--}{\cdot}  \wt k_{--} = 0\ ,\qquad  e^o_{-a}{\cdot} \wt k_{--} = 0$$
and that  such a solution $\wt k_{--}$  can be assumed to take any desired valued   at the  points of $\cV \times \cV \times \{U_o\} \times G$.  This can be checked by observing that an  explicit expression for any holomorphic  solution  $\wt k_{--}$ is actually determined in \cite[Lemma 5.3]{DS} as   
$$\wt k_{--}(x, y,  U, g) =\displaystyle  \sum_{\smallmatrix m,n \ge 0\\ m+n = 2\endsmallmatrix}   c_{mn}(x, y, g) (u^1_{-})^{m} (u^2_-)^{n} \ ,$$
 with $U = (u^i_{\pm})$. Observing that the matrix $U_o = I_2$ has entries $u^1_{-}$,  $u^2_{-} $ equal to $0$, $1$, respectively, by appropriately choosing the component $c_{0 2}(x, y , g)$, the function can take  
 any desired value  at the points  
 of $\cV \times \cV \times \{U_o\} \times G$. 
  Since $\f^o((\cV\times \cV \times \{U_o\} \times G) \cap \wh S) \subset (\cV \times \cV\times \{U_o\} \times G) \cap \wh S$,   such a residual degree of freedom  for the $\wt h_{--}$ can be used to make  the restriction $h_{--} |_{\wh S \cap \cT^{(x,U_o, g)}}$  identically zero. 
\par
 We combine the  solutions along the leaves $\cT^{(x,y, U_o, g)}$, $(x, y, U_o,  g) \in \wh S \times G$,   into a global solution
 of   \eqref{original} on $\cV \times \cV \times \SL_2(\bC) \times G$ and   restrict  such a globally defined  
 $\gg$-valued map $h_{--}$   to the submanifold $\cV \times \cV \times \SL_2(\bC) \times \{e\} \simeq \cV \times \cV \times \SL_2(\bC)$.  Then, along each  integral leaf in $ \cV \times \cV \times \SL_2(\bC)$ of the distribution spanned by   $H_{\a}$ and $e_{-a}$,  we may  consider  a new $\gg$-valued function $\psi'$ satisfying the differential problem 
$$ H_0{\cdot}  \psi' = 0 \ ,\qquad H_{--}{\cdot}  \psi' =  h_{--} \  , \qquad e_{-a}{\cdot} \psi' = 0\ .
$$
The same argument as before shows the existence of such a solution. Moreover, the residual degree of 
freedom in the choice of the solution may be used to set it to  $0$ at  each point of the form $(x, y, U_o)$. Combining  the solutions along  all the considered integral  leaves, we get a global solution $\psi'$ such that $\psi'(x, y, U_o) = 0$ for all $(x, y) \in \cV \times \cV$ and satisfying the differential problem 
\beq   \label{expression0} H_0{\cdot}  \psi' = 0 \ ,\quad H_{--}{\cdot}  (H_{--}{\cdot} \psi' )=  - H_{--}{\cdot} A^o_{--} - [H_{--}{\cdot} \psi', A^o_{--}]\  , \quad e_{-a}{\cdot} \psi' = 0\ .\eeq
Let us now consider the  new gauge $\f$, obtained  by applying to the  analytic gauge $\f^o$  the gauge transformation $g_{(x,y)} = e^{\psi'_{(x,y)}}$, $(x, y) \in \cV \times \cV$. Recalling that $A^o_0  = A^o_{-a}  = 0 $, we see that 
 the potential  in  the  new gauge $\f$ has components $A_0 $, $A_{--} $,  $A_{-a}$ given by   
\begin{multline}  \label{expression} A_0 = e^{-\ad_{\psi'}}(A^o_0 + H_0{\cdot} \psi') = 0\ ,\qquad A_{-a} = e^{-\ad_{\psi'}}(A^o_{-a} + e_{-a}{\cdot} \psi') = 0\ ,\\
A_{--} = e^{-\ad_{\psi'}}(A^o_{--} + H_{--}{\cdot} \psi')\ .
\end{multline}
From this and \eqref{expression0}, it follows that   
$$ H_{--}{\cdot} A_{--} =  e^{ad_{\psi'}}\left( [H_{--}{\cdot }\psi',  A^o_{--} ] + H_{--}{\cdot} A^o_{--}  +   H_{--}{\cdot}(H_{--}{\cdot} \psi' )\right) = 0$$
and, with a similar computation,  that $ H_0{\cdot} A_{--} = - 2 A_{--}$ and $e_{-a}{\cdot} A_{--} = 0$. 
In other words, the new prepotential $A_{--}$ satisfies  all three conditions in the first line of  \eqref{4.13}. To check the  last two normalising conditions, 
 we  recall at first  that  in any central frame, the corresponding potential $\wt A$ satisfies
 $$H_0{\cdot} \wt A_{-a}  =  -\wt A_{-a} \ ,\qquad H_{++}{\cdot} \wt A_{-a} = \wt A_{+a}\ ,\qquad H_{--}{\cdot} \wt A_{-a} = 0\ .$$ 
 Hence,   if we differentiate the identity $\psi^o|_S = 0$ in  directions  tangent to $S$,     using 
   \eqref{problemino},  \eqref{4.13} and   commutation relations we get 
\beq\label{derivpsi}
e_{+a}{\cdot} \psi^o|_S = H_{++}{\cdot} \psi^o|_S  = H_{--}{\cdot} \psi^o|_S +  \l^a \wt A_{-a}|_S = 0\ .
 \eeq
Since  $e^{\psi^o}|_S = \Id$,  we conclude that  
$A^o_{--}|_S = \wt A_{--}|_S + H_{--}{\cdot} \psi^o|_S = -   \l^a \wt A_{-a}|_S$. 
Now, from \eqref{expression0},  the property that $\wh S \subset S$ and 
$\psi'|_{\wh S} = H_{--}{\cdot} \psi'|_{\wh S} = 0$, the second line  in \eqref{4.13} follows immediately. 
\end{pf}
\par
The  $\gg$-bridge $\psi$, the  corresponding analytic gauge $\f$ and the associated prepotential $A_{--}$  established  by Lemma \ref{lemmabridge} are called {\it normalised at  $(x_o, x_o, U_o)$}. \par

\begin{rem} \label{remark47} The proof of the previous theorem shows that {\it any prepotential} $A_{--}$  (not just the normalised ones) is such that $X_-{\cdot} A_{--} = 0$ for any horizontal vector field in the distribution $\cD_-$. On the other hand, we also have that: 
\begin{itemize}
\item[a)]  $F(V, \cdot) = 0$ for any vertical  vector field $V$ of  $\HCM$, i.e. for any $V$ which is  tangent to the vertical fibers $\{(x,y)\} \times \SL_2(\bC)$ of $\HCM$.
\item[b)] The vertical anti-holomorphic distribution  of  $\HCM$ is spanned by right-invariant vector fields along the fibers $\{x\} \times \SL_2(\bC)$, which therefore commute with the left-invariant holomorphic vector field  $H_{--}$.
\item[c)] The holomorphic potential $A$ in an analytic gauge vanishes identically along the anti-holomorphic vector fields of $\HCM$. This is  due to  Proposition \ref{prop41} and the fact that the analytic gauges  are holomorphic with respect to the complex structures of  the extended  $P'$ over $(\HCM, \bI^\bC)$.   
\end{itemize}
From (a), (b), (c) and the explicit expression of $F$ in terms of a potential it follows that,  for any  anti-holomorphic vertical vector field  $V^{01}$ of $\HCM$, we have  $V^{01}{\cdot} A_{--}  = 0$. This and the   above property  $X_-{\cdot} A_{--} = 0$  prove that {\it for any prepotential $A_{--}$ defined on some  open set $\cW$ of the complexified harmonic space  $(\HCM, \bI^\bC)$,   the restriction  $A_{--}|_{\HM \cap \cW}$ is also holomorphic with respect to the complex structure $\bI$ of the harmonic space $\HM$}. \end{rem}
\par
\section{Existence,  uniqueness and compactness theorems}
\setcounter{equation}{0}
\subsection{Existence and uniqueness of an instanton  with a given prepotential}
 \begin{theo} \label{theorem51}   Let  $\cV \subset M$ be open and simply connected. 
For any    map $A_{--}: \cV \times \SL_2(\bC) \subset \HM \to \gg$, which  is holomorphic 
(i.e. with $X_-{\cdot} A_{--} = 0$ for any $X_- \in \cD_{-}$ and holomorphic in the complex coordinates of $\SL_2(\bC)$) and satisfies
\beq \label{cond-prepotential} H_0{\cdot} A_{--} = - 2 A_{--}\ ,
\eeq
there is a unique instanton $(E|_{\cV}, D)$ on $\cV$   and  an analytic gauge   $ \f\colon\!\!  P|_{\cV {\times} \SL_2(\bC)} \!\! \to$ $(\cV {\times} \SL_2(\bC)) {\times} G$, for which  $A_{--}$ is  the  prepotential in that gauge. 
\end{theo}
  \begin{pf} In order to prove the existence of an associated instanton, we  proceed as in  the proof of  \cite[Thm. 4]{ACD} and we  consider an orbit $\{x\} \times \cO \subset \cV \times \SL_2(\bC)$ of the Borel subgroup $B \subset \SL_2(\bC)$, generated by   $\langle H^o_{--}, H^o_0  \rangle \subset \sl_2(\bC)$, i.e. 
  $$B = \left\{\  \left(\begin{array}{cc} \z & 0 \\ z & \z^{-1} \end{array}\right)\ , \ (\z, z) \in \bC^* \times \bC\ \right\} \subset \SL_2(\bC)\ . $$
Along such an orbit, we may  consider the unique  $\gg$-valued connection for the $G$-bundle $\{x\} \times \cO \times G$, whose  $(1,0)$-potential   is defined by  $ A(H_{--}) = A_{--}$ and $A(H_0) \= 0$. 
 Due to  \eqref{cond-prepotential} and the fact that $[H_0, H_{--}] = - 2 H_{--}$,  such a connection has  zero curvature.  Hence there exists a new holomorphic  gauge   $\f : \{x\} \times \cO \times G \to \{x\} \times \cO \times G$  which 
  fixes all points of   $\{x\} \times \cO \times \{e\}$   and transforms $A$  into the  identically vanishing potential. By  \eqref{changepot}, 
 this is tantamount to saying that  the associated gauge transformation  
 $g_{(x, U)} $, $(x, U)  \in \{x\} \times \cO$,  is a solution to the differential problem
 \beq\label{eq1}   H_0{\cdot} g= 0\ ,\qquad H_{--}{\cdot} g +  A_{--} = 0\ , \qquad g|_{(x, U)} = e\ .\eeq 
 Since the space of all $B$-orbits $\{x\}{\times}\cO$ in $\cV{\times}\SL_2(\bC)$ is diffeomorphic to $\cV \times \bC P^1$ and is therefore simply connected, all these new gauges combine into a globally defined gauge   
 \beq \label{bri} \f  : \cV   \times \SL_2(\bC)\times G \to \cV \times \SL_2(\bC) \times G\ ,\eeq
which  maps each $B{\times}G$-orbit into  itself  and satisfies  \eqref{eq1} at all points. We may now consider   the $\gg$-valued map on $\cV \times \SL_2(\bC)$
  $$A_{++}: \cV\times \SL_2(\bC) \to \gg\ ,\qquad A_{++} \=   - (H_{++}{\cdot} g) g^{-1}\ .$$
Combining this function with  the  $\gg$-valued maps $A_{--}$ and $A_0 = 0$, we may construct the
 $G$-invariant  vector fields $H^h_0$, $H^h_{\pm\pm}$ on $\cV \times \SL_2(\bC) \times G$, which on the submanifold  
 $\cV \times \SL_2(\bC) \times\{e\}$ are 
  \beq
 \begin{split}
 \label{Hhoriz} & H^h_0|_{\cV\times \SL_2(\bC) \times \{e\}} \= H_0|_{\cV \times \SL_2(\bC) \times \{e\}}\ ,\\
& H^h_{\pm\pm}|_{\cV\times \SL_2(\bC)  \times \{e\}} \= H_{\pm \pm}|_{\cV\times \SL_2(\bC)  \times \{e\}} + A_{\pm\pm}\ \ .
\end{split}
\eeq
 By real analyticity, there is an open $\cU\,{\subset}\,\HCM$   with   $\cU \cap \HM = \cV{\times}\SL_2(\bC) $, where these  vector fields extend as holomorphic  fields on $\cU{\times}G$. Moreover, by the construction of $A_{--}$ and of the map \eqref{bri}, along each fiber  $\{x\} \times \SL_2(\bC)\times G$, the functions $A_0 \= 0$, $A_{\pm\pm}$ can be considered as  the three components of the holomorphic potential of a   connection for the $G$-bundle $\pi^{(x)}\colon \{x\} {\times} \SL_2(\bC){\times} G \to \{x\} {\times}\SL_2(\bC)$,  which is   transformed  by   the gauge  $\f$   into the trivial potential. This means that 
 the  associated covariant derivative has   identically vanishing  curvature, i.e.  $F(H_\a, H_\b) \equiv 0$ for all  $\a, \b \in \{0, ++, --\}$. \par
 \smallskip
These connections on the  submanifolds $\{x\} \times  G$ of the open set  $\cU \times \SL_2(\bC) \times G$  can be considered as   restrictions of 
a $G$-connection on  the (trivial) bundle $p: \cU \times G \to \cU$, with  associated holomorphic 
  potential $A: \cU \to  T^* \cU \otimes \gg$   satisfying 
   \begin{align}
 \nonumber & A(H_\a) \= A_\a\ ,\qquad &&\a \in  \{0, ++, --\}\ ,\\
 & A(X_-) \= 0\ ,\qquad &&\text{for any }\ X_-   \in \cD_{-}\ ,\\
 \nonumber &  A(e_{+a}) \= - e_{-a}{\cdot}  A_{++}\ \ &&\text{for any frame field}\ (e_{+a}, e_{-b})\  \text{generating}\ \cD_+ \oplus \cD_-\\
 \nonumber & \hskip 1 cm &&\text{and projecting onto the complex frames $(e_{\pm a})$ of} \ M\\
 \nonumber & \hskip 1 cm &&\text{described in \S  \ref{adapted} (for their existence, see below)}.
 \end{align}
We now need to show that  the curvature of the associated covariant  derivative $D'$  satisfies the following equalities for  each  $X_\pm, Y_\pm \in \cD_\pm$:
\beq
 \begin{split}
 \label{55} &F(X_+, Y_+) = 0 = F(X_-, Y_-) \ ,\\
  & F(H_\a, X_+) = 0 = F(H_\a, X_-)\ , \quad\ \a \in \{0, ++, --\}\ .
\end{split}
\eeq
 This would conclude the proof. Indeed,   by Proposition \ref{prop28} and Lemma \ref{char},  it would imply that  the corresponding gauge field $(E'|_{\cU} , D')$ is  the extension to   $\cU \subset \HCM$ of the  lift  of an instanton  $(E|_{\cV} , D)$ on $\cV$ and with    $A_{--}$  as  prepotential. Further, by the above definition of the potential, the trivial gauge $\f = \Id_{\cV \times \cV \times  \SL_2(\bC) \times G}$ would be  an analytic gauge for such an instanton.  Hence, by Theorem \ref{thm53}, any  instanton having  prepotential  $A_{--}$  in such a   gauge would  necessarily  coincide  with  $(E|_{\cV} , D)$. 
   \par
In order to prove  \eqref{55},  for each given  point $(x_o, y_o, U_o)  \in \cU  \subset \HCM = M \times M \times \SL_2(\bC)$  we  select  a (locally defined)  collection of $\bI^{(\bC)}$-holomorphic  vector fields $(e_{+a}, e_{-a})$, which  generate   the distribution  $\cD_+ \oplus \cD_-$ and project  to  the  complex vectors in $T^\bC M$ determined by  adapted frames of $M$ and  described     in \S \ref{adapted}. A collection  of vector fields $e_{\pm a}$ of this kind  can be  constructed through the  following three step procedure: 
\begin{itemize} [itemsep=2pt, leftmargin=18pt]
\item[1)] Pick a local section $\s: \cU' \to O_g(\cV \times \cV, J_\a)$ of the bundle of adapted frames over $\cV \times \cV$.
\item[2)] Consider the  restrictions to the  section $\s(\cU')$  of the  canonical horizontal vector fields $\n_{e_{\pm a}} \= e_{\pm a} + \G_{\pm a}{}_I^J$  determined by the Levi-Civita connection (here, the  $\G_{\pm a}{}_I^J$ are the  Christoffel symbols).
\item[3)] Take the projections $(e_{\pm a})$ of the vector fields  $\n_{e_{\pm a}} \= e_{\pm a} + \G_{\pm a}{}_I^J$ onto the underlying open set $\cU'$.
\end{itemize}
 Since the  Levi-Civita connection has vanishing torsion,  the vector fields $e_{\pm a}$ constructed  in this way are such that 
\beq\label{2} [e_{+a}, e_{+b}] = 0 = [e_{-a}, e_{-a}]\ .\eeq
Further,    just by looking at   the standard    action of the   $H^o_\a \in \sl_2(\bC)$  on  the    elements  $(h^o_i \otimes e^o_a) \subset \bC^2 \otimes \bC^{2n} \simeq (\bH^n)^\bC$, 
the   actions of the  vector fields $H_\a$ of $\HCM = M \times M \times \SL_2(\bC)$ on the   $e_{\pm a}$ are  
\beq\label{1} [H_0, e_{\pm a}] = \pm e_{\pm a}\ ,\qquad [H_{\pm\pm}, e_{\pm a}] = 0\ ,\qquad [H_{\pm\pm}, e_{\mp a}] = e_{\pm a}\ .\eeq
 Hence, by \eqref{Fcurvature} and the assumption  $A(H_0) = A(e_{-a}) = 0$, 
\beq \label{5.7} F(H_0, e_{-a}) = F(e_{-a}, e_{-b})  = 0\ .\eeq
 From  condition 
$A(e_{+a}) = -e_{-a}{\cdot} A(H_{++})$,  we  also have   
\beq  \label{5.8}
\begin{split}  F(H_{++}, e_{- a}) &=   - e_{-a}{\cdot} A(H_{++}) - A(e_{+a}) = 0 \ , \\
F(H_{--}, e_{-a}) & =- e_{-a}{\cdot} A(H_{--})  = - e_{-a}{\cdot} A_{--}  \overset{\eqref{cond-prepotential}}= 0 \ .
\end{split}
\eeq
 Finally, from    \eqref{5.7},  \eqref{5.8}, the property $F(H_\a, H_\b) = 0$ and    Bianchi identities  
\beq \label{bianchi} \sum_{\smallmatrix \text{cyclic permutations}\\ \text{of}\ (1,2,3) \endsmallmatrix} D'_{X_{i_1}} F(X_{i_2}, X_{i_3}) +  F(X_{i_1},[X_{i_2}, X_{i_3}]) = 0\eeq
 with   $X_1, X_2, X_3 $ equal to   the triple   $H_{++}, H_{--}$,   $ e_{-a}$   or to the triple $H_{++}, H_0$,   $e_{-a}$,  we  get that 
$ F(H_{--}, e_{+b}) = 0 = F(H_{0}, e_{+b})$. All this shows that   
\beq
\begin{split}
& F(e_{-a}, e_{-b})  = 0\ ,\\
\label{511} &  F(H_\a, e_{-a}) = 0\ ,\quad  \a \in \{0, ++, --\}\ ,\\
& F(H_0, e_{+a}) = F(H_{--}, e_{+a}) = 0\ .\\
\end{split}
\eeq
Now, in order to conclude the proof we  need the following 
\begin{lem} \label{lemma5.2}  The components $F(H_{++}, e_{+a})$,  $F(e_{+a}, e_{+b})$, $1 \leq a, b \leq 2n$,  of the above defined covariant derivative $D'$   are identically vanishing. 
\end{lem}
\begin{pflemma} We recall  that, by real analyticity,   $F(H_{++}, e_{+a})$,  $F(e_{+a}, e_{+b})$ 
are extended as holomorphic functions on an $\SL_2(\bC)$-invariant neighbourhood $\cU\times \{e\} \subset \HCM \times \{e\}$ of $\cV \times \SL_2(\bC) \times\{e\}$.
Note also that, due to \eqref{511}   and Bianchi identities amongst the vector fields $H_0, H_{++}$,   $ e_{+a}$  or $H_{++}, H_{--}$,   $e_{+a}$, respectively,    we have that 
\beq
\begin{split}
\label{514}  &D'_{H_0} F(H_{++}, e_{+a}) = 3 F(H_{++}, e_{+a})\ ,\\
& D'_{H_{--}} F(H_{++}, e_{+a}) = 0\ .
\end{split}
\eeq
If we consider the (holomorphic extensions of the) $G$-invariant  vector fields $H^h_0$, $H^h_{\pm\pm}$ on $\cV \times \SL_2(\bC) \times G$  defined by \eqref{Hhoriz} and the  $G$-equivariant function $$\cF_{+++a}:  \cU \times G \to \gg\ ,\quad \cF_{+++a}|_{ \cU\times \{e\}} \= F(H_{++}, e_{+a})\ ,$$
 we see that     \eqref{514} is equivalent to the  system of equations for the  $\cF_{+++a}$
\beq
\begin{split}
\label{515}  &H^h_0{\cdot} \cF_{+++a} = 3 \cF_{+++a}\ ,\\
&H^h_{--} {\cdot} \cF_{+++a} = 0\ .
\end{split}
\eeq
By \cite[Lemma 5.3]{DS} (or, more precisely, by its analogue involving  the vector field $H_{--}$ in place  of $H_{++}$) the restriction of \eqref{515} to each orbit $\cO \subset  \cU \times G$ of the $\SL_2(\bC)$-action  generated   by the  $H^h_\a$, 
admits  exactly  one solution, namely  the identically zero function.  It follows that  $\cF_{+++a}  = F(H_{++}, e_{+a}) \equiv 0$ on $\cV \times \SL_2(\bC) \times G$.\par
Let us now focus on the components $F(e_{+a}, e_{+b})$. By  \eqref{511} and Bianchi identities \eqref{bianchi}  amongst the  vector fields   $H_0, e_{+a}$,   $ e_{+b}$  and $H_{--}, e_{+a}$,   $ e_{+b}$,  we have that 
\beq
\begin{split}
\label{516}  &D'_{H_0} F(e_{+a}, e_{+b}) = 2 F(H_{++}, e_{+a})\ ,\\
& D'_{H_{--}} F(e_{+a}, e_{+b}) =  F(e_{+a}, e_{-b}) - F(e_{+b}, e_{-a})  = 0\ ,
\end{split}
\eeq
where the last equality is a consequence of the fact that $A_{-a}  = 0$ and of 
\beq \label{curv-comp} F(e_{+a}, e_{-b}) =  e_{-b}{\cdot}\left(e_{-a}{\cdot} A_{++}\right) \overset{[e_{-a}, e_{-b}] = 0} =  e_{-a}{\cdot}\left(e_{-b}{\cdot} A_{++}\right) = F(e_{+b}, e_{-a}) \ .\eeq
By the same argument as before,  the unique $G$-equivariant extension  of the $\gg$-valued function $\cF_{+ + ab} \= F(e_{+a}, e_{+b}) $ is solution to the differential problem 
$H^h_0{\cdot}  \cF_{++ab} = 2 \cF_{++ab}$ and $H^h_{--} {\cdot} \cF_{++ab} = 0$. As above, by  \cite[Lemma 5.3]{DS}, we get that  $\cF_{++ab} = F(e_{+a}, e_{+b})$ vanishes identically.
 \end{pflemma}
From \eqref{511} and  Lemma \ref{lemma5.2}, all conditions in  \eqref{55}  are satisfied.
 \end{pf}
\par
\subsection{Bounds for normalised prepotentials}
Let   $ \cV \times \SL_2(\bC)$   be an $\SL_2(\bC)$-invariant open subset of  $\HM$ 
with $\cV \subset M$ relatively compact and simply connected and such that $ \cV \times \SL_2(\bC)$ is a domain for both an  exp-central gauge $\wt \f$ and  an  analytic gauge $\f$ for an instanton $(E, D)$ which is normalised around $(x_o, x_o, U_o \= I_2) {\in}\cV {\times}\cV{\times} \SL_2(\bC)$, with $x_o \in \cV$. 
Let us also denote by 
 $A_{--}\colon  \cV \times \SL_2(\bC) {\to} \gg$ the  corresponding normalised  prepotential. 
\par
\begin{theo} \label{prelim-1} There exists a relatively compact  simply connected neighbourhood $\cV' \subset \cV$ of $x_o$, such that 
for each compact subset $K \subset \overline{\cV'} \times \SL_2(\bC)$,  there is a constant  $c_{K, \overline \cV} > 0$, depending just on  the compact sets $K$ and $\overline{\cV}$,  
 such that
  \begin{align} 
\label{6.14}  &\| A_{--}\|_{\cC^0(K, \gg )} \leq c_{K, \overline \cV}\| F\|^{(\wt \f)}_{\cC^0(\overline \cV, \gg)}\ .
\end{align}
\end{theo}
\begin{pf} 
Let $\wh S \subset \cV \times \cV \times \SL_2(\bC)$ be the complex submanifold passing through 
$(x_o, x_o, I_2)$ and defined in \S \ref{normalisation}, and   for each $(x, y, U) \in \wh S$, let   $\cS^{(x, y,U)}$  be the unique integral leaf through $(x,y, U)$ of the complex distribution generated by the  holomorphic vector fields $e_{-a}$ and  $H_0$ of $\HCM$.   Let also $\cV'$ be a relatively compact connected open subset of $\cV$ which contains $(x_o, x_o, I_2)$ and   such that the $\SL_2(\bC)$-invariant  set $\overline{\cV'}\times \overline{\cV'} \times \SL_2(\bC)$ is included in 
\beq \label{connectness} \overline{\cV'}\times \overline{\cV'} \times \SL_2(\bC) \subset \bigcup_{(x,y, U) \in  \wh S} \cS^{(x, y,U)} \cap ( \cV \times \cV \times \SL_2(\bC)) \ .\eeq
The existence of such  a $\cV'$ is guaranteed by the fact that the family of integral leaves $\cS^{(x, y,U)}$ is $\SL_2(\bC)$-invariant.  The condition \eqref{connectness} is chosen to ensure that any point of $\overline{\cV'}\times \overline{\cV'} \times \SL_2(\bC) \subset \HCM$ lies in 
some (connected) intersection   $\cS^{(x, y,U)} \cap  (\cV \times \cV \times \SL_2(\bC))$.\par
Let $K \subset \overline{\cV'} \times \SL_2(\bC) \subset \HM (\subset \HCM)$ be compact and denote  by 
$K' \subset \wh S $  the set of points $(x,y, U) \in \wh S$ such that  $\cS^{(x,y,U)} \cap K \neq \emptyset$.   
The set $K'$  is compact. 
Indeed,  it is the intersection between $\wh S$ and the compact set of the (regular) orbits of the points of $K$  by the action on $ \cV \times \cV \times \SL_2(\bC)$ of the local group generated by the holomorphic vector fields $H_0$ and  $e_{-a}$.\par
Since $A_{--}$ satisfies  \eqref{4.13}, by  integration of  the conditions  $e_{-a}{\cdot} A_{--} = 0$ and $ H_0{\cdot} A_{--} = -2 A_{--}$ along each connected intersection  $\cS^{(x, y,U)} \cap \cV \times \cV \times \SL_2(\bC)$, 
it follows that there exists a  constant $C_{\overline \cV}> 0$, depending only on  $\overline \cV$,  such that 
$$
\| A_{--}\|_{\cC^0(K, \gg)} \leq  C_{\overline \cV} \left(\sup_{\smallmatrix(x,y, U) \in K' \subset \wh S 
\endsmallmatrix } \| A_{--}\| \right) = C_{\overline \cV}  \| A_{--}\|_{\cC^0(K', \gg)} \ .
$$
On the other hand,  by the second line  in \eqref{4.13},  
$$ \| A_{--}\|_{\cC^0(K', \gg) }  \leq C'_{K'} \| \wt A\|_{\cC^0(K', \gg)}$$
 for some constant $C'_{K'}$ depending only on the set $K'$ or, equivalently, only  on the compact set $K$. From  this and   the fact that $\wt A$ does not depend on the coordinates of $\SL_2(\bC)$,  we infer that  $\| A_{--}\|_{\cC^0( K, \gg)}  \leq  C_{\overline \cV} C'_{K'} \| \wt A\|_{\cC^0(p_1(K'), \gg)}$, where $p_1: \HM = M \times \SL_2(\bC) \to M$ is the projection onto the first factor.  Since  $p_1(K') \subset  \overline \cV$, the claim  follows from  \eqref{Uhlenbeck-central}. 
\end{pf}
\par
\subsection{The second prepotential and the curvature of instantons}
Let   $ \cV \times \SL_2(\bC)$   be an $\SL_2(\bC)$-invariant open subset of  $\HM$ 
with simply connected $\cV \subset M$, which is domain for both an  exp-central gauge $\wt \f$ and  a (not necessarily normalised) analytic gauge $\f$. Further, let 
 $A_{--}:  \cV \times \SL_2(\bC) \to \gg$  be the  prepotential of  an instanton $(E', D')$ in  the analytic gauge $\f$. 
\par
\begin{definition} The {\it second prepotential  in   the analytic gauge $\f$} is  the $\gg$-valued map 
$A_{++} \= A(H_{++}):  \cV \times \SL_2(\bC) \to \gg$,  determined by the evaluation of the holomorphic  $(1,0)$-potential $A$ of the extension on $\HCM$ of the   pair $(P', \o')$   along the vector field $H_{++}$.
\end{definition}
In what follows, to keep a  clear distinction between $A_{--}$ and $A_{++}$,  we  sometimes call  $A_{--}$  
the  {\it first} prepotential. Either function yields a complete local description of instantons on hk manifolds,
since either one is completely determined by the other.  However, in our framework, some features of the
two descriptions are complementary:
\begin{itemize}
\item[1)] $A_{--}$ is    holomorphic on $(\HM, \bI)$  and is a solution to the simple first order linear equation  
$H_0{\cdot} A_{--} = -2 A_{--}$.  If no normalisation is taken, there are no further restrictions.   However,  there is  no  direct  way   to compute the curvature tensor from      $A_{--}$. 
\item[2)]  $A_{++}$ satisfies a (set of) second order nonlinear equations,  but  in terms of it the curvature  
is  given  by the simple formula,  
\beq\label{5.7bis} \cF^\f(e_{+a}, e_{-b}) = e_{-a}{\cdot}(e_{-b}{\cdot} A_{++})\ ,\eeq 
i.e. {\it the  non-trivial components of the curvature   are  just the   second order derivatives of $A_{++}$  along the anti-holomorphic  directions $e_{-a}$}. 
\end{itemize}
\par
The above-mentioned nonlinear equation for $A_{++}$ has been used in various contexts in the physics literature,
where it is known as the Leznov equation (see e.g.\ \cite{Lez,DL2, Si,DL1,DO,De}).
The  next lemma provides useful   relations between the  $\cC^k$-norms of the two types of prepotentials. \par
\begin{prop} \label{prelim} Let  $A_{--}:  \cV \times \SL_2(\bC) \to \gg$ be the (first) prepotential for an instanton and $A_{++}$ the corresponding second prepotential.  Then: 
\begin{itemize}[itemsep=2pt, leftmargin=18pt]
\item[1)]  
$A_{--}$ is the  unique   solution to the differential problem for
the unknown  $B_{--}$,
\beq \label{characterising}  H_{++}{\cdot} B_{--} = H_{--}{\cdot} A_{++} -[A_{++}, B_{--}]\ ,\qquad H_0{\cdot} B_{--} = - 2 B_{--}\ .
\eeq
A similar claim holds for $A_{++}$, provided appropriate sign changes are made.
\item[2)] For each  $k \geq 1$, there exist
 constants $M_{k,}, M_{k}{}' > 0$, depending only on $k$,   such that for $x \in \cV$, 
\begin{align}
\label{6.5*}  &\| A_{--}|_{\{x\} \times \SU_2}\|_{\cC^k(\SU_2, \gg)} \leq M_{k} \| A_{++}|_{\{x\} \times \SU_2}\|_{\cC^k(\SU_2, \gg)}\ , \\
\label{6.6**} &\| A_{++}|_{\{x\} \times \SU_2}\|_{\cC^k(\SU_2, \gg)} \leq M_{k}{}' \| A_{--}|_{\{x\} \times \SU_2}\|_{\cC^k(\SU_2, \gg)}\ .
\end{align}
\end{itemize}
\end{prop}
\begin{pf} (1) Let
$(E'|_{\cV \times \SL_2(\bC)}, D')$  be the lifted instanton,  for which $A_{\pm\pm}$,  $A_0 (= 0)$ are components of the potential in 
some analytic gauge  $\f$.  The expression in terms of curvature components of  the identity $\cF^\f(H_{++}, H_{--}) = 0$ 
shows that $A_{--}$  solves
 \eqref{characterising}. 
This solution is unique. Indeed,   if there were another solution to \eqref{characterising}, say $A'_{--}$,  
the $\gg$-valued map $\wt B_{--} \=    \Ad_{\Psi} (A'_{--} - A_{--}) $ 
(here $\Psi$ is the gauge transformation from $\f$ to an exp-central gauge)
would be  a solution to the differential problem 
$H_{++}{\cdot}  \wt B_{--} = 0$, $H_{0}{\cdot}  \wt B_{--} = - 2 \wt B_{--}$. 
This fact and   \cite[Lemma 5.3]{DS} would  imply that
$\wt B_{--} = 0$. Interchanging the signs, the 
corresponding   claim  for $A_{++}$ follows.
 \\[5pt]
(2) Consider the real basis for $\su_2 \subset \sl_2(\bC)$  given by the elements
\beq 
\begin{split}   G^o_1 \=  H^o_{++} - H^o_{--} = \begin{pmatrix} 0 & 1\\ -1 & 0 \end{pmatrix}&\ ,\ \   G^o_2 \=  iH^o_{++} + i H^o_{--}  \= \begin{pmatrix} 0 & i\\ i & 0 \end{pmatrix}\ ,\\
G^o_0 &= i H^o_0 = \begin{pmatrix} i  & 0\\ 0 & - i \end{pmatrix}\ ,
\end{split}
\eeq
and let  $G_0 = i H_0$, $G_1 = H_{++} - H_{--}$, $G_2 = i H_{++}  + i H_{--}$
be the corresponding vector fields on $\HM = M \times \SL_2(\bC)$.  Note that, for   each    $(x, U) \in M\times \SU_2$,  the {\it real} vectors $G_\a|_{(x, U)}$, $1 \leq \a \leq 3$,   give a frame for the tangent space of the {\it totally real}  submanifold $\{x\}\times \SU_2$  of $\{x\} \times \SL_2(\bC)$. \par
Let
$(E'|_{\cV \times \SL_2(\bC)}, D')$  be the lifted instanton,  for which $A_{\pm\pm}$,  $A_0$ are  components of the potential in 
some analytic gauge  $\f$. Since $A_0 = 0$ 
 and the curvature $F'$ of  $(E'|_{\cV \times \SL_2(\bC)}, D')$ vanishes identically along any vector field that is   tangent to the  fibres of $\HM$, at each   $(x, U) \in \{x\}\times  \SU_2$  we have  that 
  \beq \label{6.5}
 \begin{split}
 & 0 = i F'(H_0, H_{--}) =   i H_0{\cdot}A_{--} + 2 i A_{--} = G_0{\cdot} A_{--} +2  i A_{--}\ ,\\
&  0 =  2 F'(H_{++}, H_{--}) = 2 H_{++}{\cdot} A_{--}   - 2 H_{--}{\cdot} A_{++}  + 2  [A_{++}, A_{--}]= \\
& \hskip 2 cm =  (G_1 - i G_2){\cdot} A_{--} + (G_1 + i G_2){\cdot} A_{++} + 2 [A_{++}, A_{--}]  \ .
 \end{split}
 \eeq
Hence,  the restrictions
$V_{--} \= A_{--}|_{\{x\} \times \SU_2}, W_{++} \= A_{++}|_{\{x\} \times \SU_2}$  solve the system
    \beq\label{6.6}
 \begin{split}
 &    G_0{\cdot} V_{--}  + 2 i V_{--} = 0\ ,\\
 &    ( G_1-iG_2){\cdot} V_{--}  +  2\ad_{W_{++}} (V_{--})  =  - (G_1 + i G_2){\cdot} W_{++}\ .\\
 \end{split}
 \eeq
We claim  that if  \eqref{6.6} is considered as a system   of equations for
$V_{--}$ with coefficients determined by $W_{++}$, then it  is equivalent to a  system given by an appropriate  first order elliptic  operator $\cP$ with trivial kernel  (for definition and first properties of elliptic operators,  we refer to \cite[Appendix G]{Be}).  To check this,  consider the  subbundle $\cE$ 
 of the  bundle 
$\mathfrak P \=   \SU_2 \times (\gg\times \gg \times \gg \times \gg)  \to \SU_2$, defined by 
$$\cE = \{ (x, U; X, Y, Z, W) \in \mathfrak P\ :\quad   X = Z, \quad Y = W \}\ .$$
We consider $\cE$  equipped with the Hermitian product along the fibres  determined  by    \eqref{innerproduct}.  
Then,  define
 \begin{multline} \cP: \cC^\infty(M \times \SU_2, \cE) \longrightarrow  \cC^\infty(M \times \SU_2, \cE)\ ,\\
 \cP\left(\smallmatrix   X \\ Y \\ X \\ Y 
 \endsmallmatrix \right) \=
 \left( \smallmatrix
G_0 + 2 i I & G_1 + i G_2 + 2 \ad_{\overline{W_{++}}}  & 0 & 0\\
G_1 - i G_2  +2 \ad_{W_{++}} & G_0 - 2 i I& 0 & 0\\
0 & 0 & -  G_0 -  2 i I  &  G_1 + i G_2 + 2 \ad_{\overline{W_{++}}}  \\
0 & 0 & G_1 - i G_2  +2 & -  G_0 + 2 i I
\endsmallmatrix\right)
 {\cdot} \left(\smallmatrix   X \\ Y \\ X \\ Y 
 \endsmallmatrix \right)
 \end{multline}
 One can directly  check that $\cP$ is a first order  elliptic operator. Moreover,  
  $\ker \cP = \{0\}$. Indeed, a quadruple  $(X, Y, X, Y)$ is in $\ker \cP$ if and only if $X$ and $\overline Y$ are both solutions on $\SU_2$  to the system  of equations  for
  $\gg$-valued maps $f$ 
  \beq\label{6.6bis}
 \begin{split}
 &    H_0{\cdot} f   =  - 2  f\ ,\\
 &    H_{++}{\cdot} f  +  \ad_{W_{++}} (f)  = 0\ .\\
 \end{split}
 \eeq
 By applying the inverse   of  a gauge transformation  from a central to the analytic gauge $\f$,     the (flat) connection on $\SU_2$,   determined by the   potential $(A_0 = 0, A_{++}, A_{---})$,  is transformed into the (flat) connection
 determined by a  potential with trivial components $(A'_0 = 0, A'_{++} = 0, A'_{--} = 0)$. In particular, the  solution $f$ to  \eqref{6.6bis} is transformed into a solution $\wt f$ of the system
  \beq\label{6.6ter}
 \begin{split}
 &    H_0{\cdot} \wt f   =  - 2  \wt f\ ,\\
 &    H_{++}{\cdot}\wt  f   = 0\ .\\
 \end{split}
 \eeq
By  \cite[Lemma 5.3]{DS} (applied to each component of the matrix valued function $\wt f$),   we get that   $\wt f$ is  equal to $0$. This shows that \eqref{6.6ter} has $(X, Y, X, Y) = (0, 0, 0, 0)$ as the
unique solution.\par
\smallskip
We may now conclude the proof of (2). In fact,  it suffices  to  observe  that \eqref{6.6}  is equivalent to 
 saying that the   section of   $\cE$, given by the 
 quadruple 
$( V_{--} , $ $\overline{ V_{--} }$,$V_{--}, \overline{V_{--}})$,  is a  solution to the differential problem 
\beq \cP \left(\smallmatrix  V_{--} \\ \overline{ V_{--} }\\ V_{--} \\\overline{V_{--}}
 \endsmallmatrix \right) =  \left(\smallmatrix   - (G_1 - i G_2){\cdot} \overline{W_{++}} \\  - (G_1 +  i G_2){\cdot} W_{++} \\ - (G_1 - i G_2){\cdot} \overline{W_{++}} \\ - (G_1 +  i G_2){\cdot} W_{++} 
 \endsmallmatrix \right)
 \eeq
Since $\cP$  is elliptic 
with $\ker \cP = \{0\}$,  classical Schauder estimates  (see e.g. \cite[Appendix H]{Be}) imply that, for each  
$k \geq 1$ there are  constants
$N_k$,  $M_k > 0$ such that 
$$\|V_{--}\|_{\cC^k(\SU_2, \gg)} \leq N_k \left \| \left(\smallmatrix   - (G_1 - i G_2){\cdot} \overline{W_{++}} \\  - (G_1 +  i G_2){\cdot} W_{++} \\ - (G_1 - i G_2){\cdot} \overline{W_{++}} \\ - (G_1 +  i G_2){\cdot} W_{++} 
 \endsmallmatrix \right)\right\|_{\cC^{k-1}(\SU_2, \gg)}   \leq M_k \|W_{++}\|_{\cC^k(\SU_2, \gg)}\ .$$
This gives   \eqref{6.5*}.  The proof of \eqref{6.6**} is similar.
\end{pf}
\par
\subsection{The local compactness  theorem}
To conclude this paper, as an example of the utility of the harmonic space formulation, we present a
streamlined proof of Uhlenbeck, Nakajima and Tian's celebrated local compactness theorem for  
Yang-Mills fields in the specific case of hk instantons.

From the  classical estimates in \cite{Uh1, Nak} (see also  \cite{Ti, Zh}),  we know that for any  geodesic ball  
$B_R = B_R(x_o)$ of radius $R$ of an $m$-dimensional  Riemannian manifold $(M, g)$,  
the   $\cC^0$-norms of  curvatures  of   Yang-Mills fields   are controlled by 
their corresponding   $L^2$- or  $L^{\frac{m}{2}}$-norms and $R$.
In fact, there are   constants $\ve$,  $C, c_m > 0$  such that
\begin{align}
\label{est1} &\|F\|_{L^2(B_R)} < \ve R^{m-4} & &\text{implies } & &   \|F\|_{\cC^0( B_{R/4}, \gg)} \leq \frac{C}{R^{\frac{m}{2}}}\| F\|_{L^2(B_R)} \\
\label{est2} &\|F\|_{L^{\frac{m}{2}}(B_{R})} < c_m  & & \text{implies} & &  \|F\|_{\cC^0( B_{\frac{R}{2}}, \gg)}\leq \frac{2^m C}{R^{n}}\| F\|_{L^2(B_{\frac{R}{2}})} 
\end{align}
Combining   these estimates with   Theorem \ref{prelim-1} yields   \par
\begin{theo}[Local Compactness Theorem for instantons on hk manifolds] Let   $B_R = B_R(x_o) \subset M$ be   a  geodesic ball  in a $4n$-dimensional hk manifold $(M, g, J_\a)$, which is  included in a  relatively compact neighbourhood $\cV'$ of $x_o$ where   Theorem \ref{prelim-1} holds.  Further let   
$(E|_{B_R} = B_R \times V, D^{(k)})$ be a sequence of (trivialised) instantons,
each with the same compact structure group $G^o$  and corresponding normalised prepotential
$A_{--}^{(k)}$ on $B_R \times \SL_2(\bC)$.\par
If the curvatures are such that  $\|F^{(k)}\|_{L^2(B_R)} < \ve R^{4(n-1)}$ for all $k$, with $\ve > 0$ as in  \eqref{est1}, then there exists a subsequence  $(E|_{B_R}, D^{(k_n)})$,  whose  curvatures  $F^{(k_n)}$   
converge uniformly to the curvature of  a limit instanton $(E|_{B_R}, D^{(\infty)})$. 
The same conclusion holds if the  curvatures are such that  $\|F^{(k)}\|_{L^{2n}(B_{2R})} < c_{4n}$ with  constant $c_{4n} > 0$ as in \eqref{est2}.
\end{theo}
\begin{pf} By \eqref{est1}, \eqref{est2}, \eqref{6.14},  the  sequence of normalised holomorphic prepotentials $A^{(k)}_{--}$ is  uniformly bounded  on any compact subset  $K$ of  $B_R {\times} \SL_2(\bC)$. It follows from  Montel's Theorem  that there is a subsequence $A^{(k_n)}_{--}$   converging uniformly on compacta  to a holomorphic map $A^{(\infty)}_{--}$, which  is the prepotential of some instanton due to  Theorem \ref{theorem51}. Using \eqref{6.6**},  we may also assume that  the   second prepotentials $A^{(k_n)}_{++}$ and all their derivatives  converge uniformly on  compacta to the second prepotential 
$A^{(\infty)}_{++}$ and its derivatives corresponding to the instanton determined  by $A^{(\infty)}_{--}$.  Thus,  by \eqref{5.7bis},   the curvatures converge uniformly on compacta as well. \end{pf}

\vskip 1cm
\hbox{\parindent=0pt\parskip=0pt
\vbox{\baselineskip 9.5 pt \hsize=3.1truein
\obeylines
{\smallsmc
Chandrashekar Devchand
Max-Planck-Institut f\"ur Gravitationsphysik 
(Albert-Einstein-Institut)
Am M\"uhlenberg 1 
D-14476 Potsdam 
Germany
}\medskip
{\smallit E-mail}\/: {\smalltt devchand@math.uni-potsdam.de
}
}
\hskip  -3truemm
\vbox{\baselineskip 9.5 pt \hsize=3.7truein
\obeylines
{\smallsmc
Massimiliano Pontecorvo
Dipartimento di Matematica e Fisica
Universit\`a di Roma III
Largo San Leonardo Murialdo 1
I-00146 Roma
Italy
}\medskip
{\smallit E-mail}\/: {\smalltt max@mat.uniroma3.it}
}
}
\vskip 0.7cm
\hbox{\parindent=0pt\parskip=0pt
\vbox{\baselineskip 9.5 pt \hsize=3.1truein
\obeylines
{\smallsmc
Andrea Spiro
Scuola di Scienze e Tecnologie
Universit\`a di Camerino
Via Madonna delle Carceri 
I-62032 Camerino (Macerata)
Italy
}\medskip
{\smallit E-mail}\/: {\smalltt andrea.spiro@unicam.it
}
}
}

\begin{thebibliography}{43}




\bibitem{ACD} D.\,V. Alekseevsky, V. Cortes and C. Devchand,
\textit{Yang-Mills connections over manifolds with Grassmann structure}, J. Math. Phys. \textbf{44} (2003), 6047--6074

\bibitem{At} M.\,F. Atiyah, Geometry of Yang-Mills fields, {\it   Scuola Normale Superiore Pisa, Pisa}, 1979

\bibitem{ADHM} M.\,F.  Atiyah, N.\,J. Hitchin, V.\,G. Drinfel'd and
              Yu.\,I. Manin, {\it Construction of instantons}, Phys. Lett. A  {\bf 65} (1978), 185--187


\bibitem{ahs}
  M.\,F. Atiyah, N.\,J. Hitchin and I.\,M. Singer,
{\it Self-duality in four-dimensional {R}iemannian geometry},
 Proc. Roy. Soc. London Ser. A  {\bf 362} (1978), 425--461

\bibitem{aw}
  M.\,F. Atiyah and R.\,S. Ward,
\textit{Instantons and algebraic geometry}, 
  Comm. Math. Phys. {\bf 55} (1977), 117--124




\bibitem{bpst}
  A.\,A. Belavin, A.\,M. Polyakov, A.\,S. Schwartz and Y.\,S. Tyupkin,
{\it Pseudoparticle xolutions of the Yang-Mills equations},
  Phys. Lett. B {\bf 59} (1975), 85--87 
  
  \bibitem{Be} A. Besse, Einstein manifolds, {\it Springer, Berlin}, 1987
  
\bibitem{BLS}   
J.\,P. Bourguignon, H.\,B. Lawson and J. Simons,
{\it Stability and gap phenomena for Yang-Mills fields} 
 Proc. Nat. Acad. Sci.  U.S.A. {\bf 76} (1979), 1550--1553



\bibitem{cdfn} E. Corrigan, C. Devchand,  D.\,B. Fairlie and J. Nuyts,
{\it First-order equations for gauge fields in spaces of dimension greater than four},
 Nuclear Phys. B  {\bf 214} (1983), 452--464

\bibitem{cfgy} E. Corrigan,  D.\,B. Fairlie,  P. Goddard and R.\,G. Yates,
{\it The construction of self-dual solutions to $\SU(2)$ gauge theory},
  Comm. Math. Phys.  {\bf 58} (1978), 223--240

\bibitem{cgk} E. Corrigan, P. Goddard and A. Kent, 
{\it Some comments on the ADHM construction 
in $4k$ dimensions}, Comm. Math. Phys. {\bf 100} (1985), 1--13 

\bibitem{De}
  C. Devchand,
{\it Oxidation of self-duality to $12$ dimensions and beyond},
  Comm. Math. Phys. {\bf 329} (2014), 461--–482

\bibitem{DL1}
  C. Devchand and O. Lechtenfeld,
{ \it Extended self-dual {Y}ang-{M}ills from the {$N=2$} string},
  Nuclear Phys. B   {\bf 516} (1998), 255--272

\bibitem{DL2}
  C. Devchand and A.\,N. Leznov,
{\it B\"{a}cklund transformation for supersymmetric self-dual theories
              for semisimple gauge groups and a hierarchy of {$A_1$}
              solutions}, 
  Comm. Math. Phys. {\bf 160} (1994), 551--561
     
     
   \bibitem{DO}
  C. Devchand and V. Ogievetsky,
{\it Interacting fields of arbitrary spin and {$N>4$}
              supersymmetric self-dual {Y}ang-{M}ills equations},
 Nuclear Phys. B {\bf 481} (1996), 188--214

\bibitem{DS} C. Devchand and A. Spiro, {\it 
On pseudo-hyperk\"{a}hler prepotentials},
J. Math. Phys.  {\bf 57}  (2016), 102501--102537 
 
 \bibitem{Do}  
  S.\,K. Donaldson,
{\it  Anti self-dual Yang-Mills connections over complex algebraic surfaces and stable vector bundles},
   Proc. Lond. Math. Soc. {\bf 50} (1985), 1--26
 
\bibitem{DK} S.\,K. Donaldson and P.\,B. Kronheimer, The geometry of four-manifolds, {\it The Clarendon Press, Oxford University Press, New York}, 1990


\bibitem{DT}
  S.\,K. Donaldson and R.\,P. Thomas,
 {\it Gauge theory in higher dimensions}, 
 in  "{The geometric universe ({O}xford, 1996)}", pp. 31--47, 
 {\it Oxford Univ. Press, Oxford}, 1998.


\bibitem{DM}
  V.\,G. Drinfeld and Y.\,I. Manin,
{\it  A description of instantons},
  Comm. Math. Phys.  {\bf 63} (1978), 177--192
 


\bibitem{FU} D.\,S. Freed and K. Uhlenbeck,
   {Instantons and four-manifolds},
{\it Springer-Verlag, New York}, {1991}

\bibitem{gios}
A. Galperin, E. Ivanov, V. Ogievetsky and E. Sokatchev,
\textit{Gauge field geometry from complex
and harmonic Analyticities.
I. K\"ahler and self-dual Yang-Mills cases},
Ann.  Physics {\bf 185}  (1988), 1--21


\bibitem{GIO}
 A. Galperin, E.  Ivanov and  O. Ogievetsky, 
\textit{Harmonic space and quaternionic manifolds}, 
Ann. Physics {\bf 230}  (1994),  201--249



\bibitem{gios_book}
A.\,S. Galperin, E.\,A. Ivanov, V.\,I. Ogievetsky and E.\,S. Sokatchev,
Harmonic superspace, {\it Cambridge University Press, Cambridge}, 2004


\bibitem{HL} R. Harvey and H. Blaine Lawson Jr., 
{\it Calibrated geometries}, Acta Math.  {\bf 148} (1982),
47--157
    
\bibitem{Le} C. LeBrun,
{\it Quaternionic-K\"ahler manifolds and conformal geometry}, 
Math. Ann. {\bf 284} (1989),  353--376

\bibitem{Lez}
  A.\,N. Leznov,
{\it On the equivalence of four-dimensional self-duality equations
              to a continuous analogue of the problem of a principal chiral
              field}, 
   Teoret. Mat. Fiz.. {\bf 73} (1987), 302 -- 307 (English translation: Theoret. and Math. Phys. {\bf 73} (1987),  1233–1236) 


\bibitem{MS}  M. Mamone Capria and S.\,M. Salamon,  {\it Yang-{M}ills fields on quaternionic spaces},
  Nonlinearity {\bf 1} (1988), 517--530 
  
  \bibitem{Mo} J.\,W. Morgan, {\it An introduction to gauge theory} in ``
  Gauge theory and the topology of four-manifolds (Park
              City, UT, 1994)'', pp. 51--143,
 {\it Amer. Math. Soc., Providence, RI}, 1998

  \bibitem{Nak} H. Nakajima, {\it Compactness of the moduli space of Yang-Mills connections in higher dimensions}, 
  J. Math. Soc. Japan\
 {\bf 40}  (1988), 383--392
  

\bibitem{Ni} T. Nitta, {\it Vector bundles over quaternionic K\"ahler manifolds}, 
Tohoku Math. J. {\bf 40} (1988), 425--440

\bibitem{Pa} R. Palais,  {\it A global formulation of the Lie theory of transformation groups},
Mem. Amer. Math. Soc. No. {\bf 22}, (1957), pp. iii+123


\bibitem{polyakov}
  A.\,M. Polyakov,
{\it Quark confinement and topology of gauge groups},
  Nuclear Phys. B {\bf 120} (1977),  429--458 


\bibitem{Sa} S.\,M. Salamon, {\it Quaternionic K\"ahler manifolds}, Invent. Math. \textbf{67} (1982), 143--171

\bibitem{Si}
  W. Siegel,
{\it {$N=2\ (4)$} string theory is self-dual {$N=4$} {Y}ang-{M}ills
              theory},
  Phys. Rev. D (3)  {\bf 46} (1992),  R3235--R3238

\bibitem{Ti} G. Tian, {\it Gauge theory and calibrated geometry, I}, Ann. of Math. {\bf 151} (2000), 193--208

\bibitem{Uh1} K.\,K. Uhlenbeck, {\it Removable singularities in Yang-Mills fields}, Comm. Math. Phys. {\bf 83} (1982), 11--29.

\bibitem{Uh2} K.\,K. Uhlenbeck, {\it Connections with $L^p$ bounds on curvature}, Comm. Math. Phys. {\bf 83} (1982), 31--42

\bibitem{UY}
K.\,K. Uhlenbeck and  S.\,T. Yau, 
{\it On the existence of Hermitian-Yang-Mills connections in stable vector bundles},
Comm. Pure Appl. Math., {\bf 39} (1986), S257--S293


\bibitem{ward77}
  R.\,S. Ward,
 {\it On self-dual gauge fields},
  Phys. Lett. A {\bf 61} (1977), 81--82
  
\bibitem{ward84} R.\,S. Ward, {\it Completely solvable gauge-field equations in dimension greater
              than four}, Nuclear Phys. B, {\bf 236},
(1984),  381--396

\bibitem{We} K. Wehrheim, Uhlenbeck Compactness, {\it European Mathematical Society, Z\"urick}, 2004


\bibitem{Zh}  X. Zhang, {\it A compactness theorem for Yang-Mills connections}, Canad. Math. Bull. {\bf 47} (2004), 624--634


\end{thebibliography}
\end{document}